\documentclass[11pt, oneside]{article}   	% use "amsart" instead of "article" for AMSLaTeX format
\usepackage{geometry}                		% See geometry.pdf to learn the layout options. There are lots.
\usepackage{graphicx,color}				% Use pdf, png, jpg, or eps§ with pdflatex; use eps in DVI mode
\usepackage[bf,hang]{caption}
\usepackage{subcaption,url}								% TeX will automatically convert eps --> pdf in pdflatex		
\usepackage{amssymb, amsmath, amsthm}
\usepackage{authblk}

\newtheorem{theor}{Theorem}[section]

 \newtheorem{prop}{Proposition}[section]
 \newtheorem{rem}{Remark}[section]
\def\R{\mathbb{R}}
\def\N{\mathbb{N}}

\allowdisplaybreaks

\begin{document}

\title{A numerical study of an Heaviside function driven\\ degenerate diffusion equation}

\author[1]{C. Alberini}
\author[2]{R. Capitanelli}
\author[3]{S. Finzi Vita}
\date{}							% Activate to distplay a given date or no date

\affil[1]{{\it\footnotesize{Dipartimento SBAI, Sapienza Universit\`a di Roma, raffaela.capitanelli@uniroma1.it}}}
\affil[2]{{\it\footnotesize{Dipartimento SBAI, Sapienza Universit\`a di Roma, carlo.alberini@uniroma1.it}}}
\affil[3]{{\it\footnotesize{Dipartimento di Matematica, Sapienza Universit\`a di Roma, stefano.finzivita@uniroma1.it}}}

\maketitle

{\bf Abstract.}
We analyze a  nonlinear degenerate  parabolic problem whose diffusion coefficient is the Heaviside function of the distance of the solution itself from a given target function. We show that this model behaves as an evolutive variational inequality having the target as an obstacle: under suitable hypotheses, starting from an initial state above the target the solution evolves in time towards an asymptotic solution, eventually getting in contact with part of the target itself. 

We also study a finite difference approach to the solution of this problem, using the exact  Heaviside function or a regular approximation of it, showing the results of some numerical tests.

\vspace{5mm}



%%%%%%%%%%%%%%%%%%%%%%%%%%%%%%%%%%%%

\section{Introduction}

Aim of the paper is to study  the following problem
\begin{eqnarray}
\label{model1}
	\left\{
	\arraycolsep=1.4pt\def\arraystretch{1.5}
	\begin{array}{ll}
		u_t -H\left(u-u^c\right)\left(\Delta u + f\right) =0 \ &\text{ a.e. in } \Omega, \text{ for all }t \in (0,T),	\\
		u\left(0\right) = u^{0} 						& \text{ in } \Omega,	\\
		u=0 & \text{ on } \partial\Omega, \text{ for all }t \in (0,T),
	\end{array}\right.
\end{eqnarray}
with $T>0$, where $H$ is the extended Heaviside function such that $H(0)=0$, that is

\begin{equation}
\label{heavi}
	H(r)=\left\{
	\arraycolsep=1.4pt\def\arraystretch{1.5}
	\begin{array}{ll} 1   &\text{ for } r>0\cr
		%[0,1]   &\text{ for } r=0 \cr 
		0  &\text{ for } r\le0
			\end{array}\right. 
\end{equation}
and $ \Omega$ is a bounded domain in $\R^n$,  $n\in\N$,  with smooth boundary.

We assume that the  initial datum $u^0$,  the (independent of time) source term $f$ and the given target function (the obstacle) $u^{c}$ satisfy the following conditions
\begin{equation}
\label{obst}
	 u^0\in H^1_0\left(\Omega\right),\ f\in L^{2}\left(\Omega\right),\ u^{c}\in H^2\left(\Omega\right),\ u^c \le 0 \text{ on } \partial\Omega.
\end{equation}

We define as solution of problem (\ref{model1}) a function $u\in L^{2}\left(0,T;H^1_0\left(\Omega\right)\cap H^2\left(\Omega\right)\right)$ with $u_t\in L^{2}\left(0,T;L^2(\Omega)\right)$ which solves problem (\ref{model1}).
 
Problem (\ref{model1})  fits into the typology of degenerate parabolic problems, since the diffusion coefficient is a discontinuous function which vanishes where the solution touches the obstacle  $u^{c}. $

The use of the Heaviside function in the formulation of evolutive differential problems is useful when a discontinuous behavior can occur according to the specific values of the solution itself: such a function  acts as a dynamic switch for this behavior. It can be applied to the differential operator itself (for example the Laplacian), giving rise to nonlocal phenomena: examples of that kind can be found in models connected to epidemics spread \cite{Mar} or self-organized criticality such as the sandpile model (see e.g. \cite{Barbu}, \cite{Mosco}, \cite{AFV}, \cite{MoscoVivaldi}). In these cases the Heaviside function, calculated on the distance between the solution and an assigned critical state,   is able  to govern the spread of the problem on a global level: the initial data tend to the critical state progressing from the edge towards the interior of the domain during the time evolution, which stops when all the solution gets in contact with the threshold  critical state. 

 Other examples of this approach can be found in the literature even without the Heaviside function, for example in cases where the model changes its behavior in a discontinuous way according to the values of the solution or of its partial derivatives: in these cases the Heaviside function is replaced by the positive part function, $\sigma\to\left(\sigma\right)_{+} := \max\left\{\sigma,0\right\}$, for all $\sigma\in\R$ to describe non-reversible phenomena that can be connected to the avalanche behavior in the sand piles \cite{Preziosi} and  to damage mechanics models \cite{Akagi}.
 
In the present paper the Heaviside function acts on the diffusion coefficient, conditioning the pointwise evolution of the solution in base of its distance from a given threshold (the obstacle function). It can henceforth be considered a useful tool when some irreversible phenomena appear   in time. The application of such formulation still arises in models for the epidemic spread (see again \cite{Mar}). 

 In Section 2 we prove  under suitable conditions the equivalence of  problem (\ref{model1}) with a parabolic obstacle problem, that is with a variational inequality on the convex set of the functions above the target. As a consequence, it is possible to  characterize the asymptotic solution of  (\ref{model1}) as the solution of the corresponding stationary (elliptic) obstacle problem.

In Section 3 we discuss a direct numerical approximation of  problem (\ref {model1}) in one and two dimensions through a semi implicit finite difference scheme, either through the use of the exact Heaviside function, or of a  $ C^1$ approximation of it; an efficient variable time step strategy is also described in order to reduce the computational costs of the method. Finally, in Section 4 we present some numerical tests.

%%%%%%%%%%%%%%%%%%%%%%%%%%%%%%%%%%%%

\section{The  model}

In the present section we analyze problem (\ref{model1}), showing that, under suitable conditions,   it is equivalent to a parabolic variational inequality with obstacle $u^c$. In such a way we will be able to prove the existence of an asymptotic solution for such a problem, and to characterize it as the solution of the corresponding stationary obstacle problem.
\subsection{Equivalence with the obstacle problem}

We start by recalling the classic parabolic obstacle problem:
\begin{eqnarray}
\label{parab_obst}
	\left\{
	\arraycolsep=1.4pt\def\arraystretch{1.5}
		\begin{array}{ll}
			\displaystyle w(t)\in \mathcal{K}, \quad \int_\Omega ( w_t - \Delta w - f)( \varphi - w)\,dx \geq 0	&	\quad\forall\varphi \in\mathcal{K}, \forall t\in\left(0,T\right)	\\
			w\left(0\right) = u^0	\in 	\mathcal{K}												
		\end{array}
	\right.
\end{eqnarray}
where $\mathcal{K}$ denotes the convex set

\begin{equation}
\label{convex1}
	 \mathcal{K} = \left\{\varphi \in H^1_0\left(\Omega\right), \varphi \geq u^c \text{ in } \Omega\right\}
\end{equation}
and $u^0, u^c, f$ are the same data used in \eqref{model1}. It is well known  (see  \cite{BR}, pp. 99-102) that under the assumptions \eqref{obst}, there exists a unique solution $w=w(x,t)$ for problem  (\ref{parab_obst}), with $w\in L^2(0,T;\mathcal{K}\cap H^2(\Omega))$ and $w_t\in L^{2}\left(0,T;L^2(\Omega)\right)$ (see also \cite{CT} %pp. 113-119
and \cite{Friedman}). %, pp. 76 - 77, Chapter 1 - section 8
Moreover (\ref{parab_obst}) can be written in the equivalent formulation of a complementarity system, that holds for all $t>0:$
\begin{eqnarray}
\label{parab_cs}
	\left\{
	\arraycolsep=1.4pt\def\arraystretch{1.5}
		\begin{array}{ll}
		w(x,t) \geq u^c(x) 								&\  	\text{in } \Omega		\\
		w_t \geq \Delta w + f							&\ 	\text{a.e. in } \Omega		\\
		\displaystyle (w-u^c)(w_t -\Delta w - f) = 0			&\	\text{a.e. in }\Omega		\\
		w\left(x,0\right)=u^0 							&\	 \text{in } \Omega		\\
		w(x,t) = 0.									&\	 \text{on } \partial\Omega
		\end{array}
	\right. 
\end{eqnarray}

We point out that  a fundamental result  for this problem is given by the parabolic  version of the  Lewy - Stampacchia inequality (see e.g.  \cite{CT}, p. 113) %and \cite{LewSt})
\begin{equation}
	\label{LS} 
	f\leq  w_t - \Delta w  \leq \sup (0,-\Delta u^c-f) + f.
\end{equation}

We are able to prove the equivalence of  problems (\ref{model1}) and (\ref{parab_obst}) under the following assumptions:
\begin{description}
 	\item[H$_{1}$:] $u^{0} > u^{c}$ a.e. in $\Omega$;
	\item[H$_{2}$:] $\Delta u^{c} + f \leq 0$ a.e. in $\Omega$.
\end{description}

We point out that  under conditions \eqref{obst} and  {\bf H}$_{2}$, if $w$ solves problem   {\em(\ref{parab_obst})} (hence  {\em(\ref{parab_cs})}), then $w$ solves {\em(\ref{model1})}.   
In fact, initial and boundary conditions are the same in (\ref{model1}) and (\ref{parab_cs}). 
If $w>u^c,$  we obtain $w_t -\Delta w - f= 0$  a.e $\text{in } \Omega.$	
	If $w=u^c,$ then by \eqref{LS} and  {\bf H}$_{2}$, $w_t - \Delta u^c  \leq \sup (0,-\Delta u^c-f) + f=-\Delta  u^c$  so that $ w_t \leq 0:$  as $w \geq u^c,$ we obtain  $ w_t = 0 .$  

Then under such conditions we obtain that a solution of problem  {\em(\ref{model1})} exists.
 In the following Proposition \ref{equiv}, we prove that under the further condition  {\bf H}$_{1}$, problems  {\em(\ref{parab_obst})} and {\em(\ref{model1})} are equivalent.

\begin{prop}\label{equiv} 
Assume that $u$ solves problem {\em(\ref{model1})}, and that  conditions \eqref{obst}, {\bf H}$_{1}$ and {\bf H}$_{2}$ hold, then $u$ coincides for any time with the unique solution $w$ of {\em(\ref{parab_obst})}, and hence of {\em(\ref{parab_cs})}.
\end{prop}

\noindent{\bf Proof.}
Initial and boundary conditions are the same in (\ref{model1}) and (\ref{parab_cs}). 

Let us now prove that if  $u(x,t)$ is a solution of \eqref{model1} then necessarily $u(.,t)\geq u^c(.)$ in $\Omega$ for any time $t$ (that is $u(t)\in \mathcal{K}$). It is true for $t=0$ thanks to {\bf H}$_{1}$. Suppose that for a given  $x\in \Omega$ there exists a first time $t^*$ such that $u(x,t^*)=u^c(x)$ from above. Then, from (\ref{model1})  $u_t(x,t^*)=0,$ so that $u(x,t)\geq u^c(x)$ for any $t>t^*.$
Then the first inequality of (\ref{parab_cs}) is satisfied by $u$.

The  equation in the third line of \eqref{parab_cs} is trivially satisfied where $u(x,t)= u^c(x)$; where $u(x,t)>u^c(x)$ then from (\ref{model1})  $u_t -\Delta u - f =0$, so it is always true. In particular we see that the detachment set of $u$ is a priori a subset of the detachment set of problem (\ref{parab_cs}). 

Concerning the second inequality of (\ref{parab_cs}), we have already seen that it is satisfied (with the equal sign) when $u> u^c$. But when $u(x,t)=u^c(x)$,  (\ref{model1}) and  assumption {\bf H}$_{2}$ imply that
$$ u_t -\Delta u - f = -\Delta u^c - f \geq 0 .$$

Then $u$ coincides with the solution of (\ref{parab_cs}); this also proves its uniqueness.\qed

 \begin{rem}\label{u0}
We point out that  the hypothesis {\bf H}$_{1}$ is crucial in order to guarantee the 
equivalence of problem {\em(\ref{model1})} and {\em(\ref{parab_obst})}. 
In fact,  if $u^0=u^c$ in some region $D\subset\Omega,$ then $u(t)=u^c$ in $D$ for any $t>0$, differently from what would happen for the solution $w(t)$ of {\em(\ref{parab_obst})}. Then  the entire evolution of the solution and hence the asymptotic solution of the problem change: 
see the next subsection discussion and an example (Test 3) of Section 4. 
\end{rem} 

\begin{rem}\label{delta_uc}
Assumption {\bf H}$_{2}$ is  a   natural condition for  the contact region  $C(t)=\{x\in\Omega: u(x,t)=u^c(x)\}$ (the place where it is used inside the proof), at least when the obstacle is sufficiently smooth. If not satisfied by the data, $C(t)$ remains empty for any time, and the two problems {\em(\ref{model1})} and {\em(\ref{parab_obst})} are trivially equivalents.  On the other side, it is possible to verify (for example,  in the case with no source term, $f=0$)  that  the contact is possible from above only at regions of $\Omega$ where the obstacle is superharmonic ($-\Delta u^c \geq 0$); in order for the solution to reach regions of the obstacle where it is subharmonic one needs a sufficiently negative source term $f$  to balance the positivity of $\Delta u^c$. Then the assumption of {\bf H}$_{2}$ in all of $\Omega$ is  too restrictive, but we left it here in this form, since the contact region itself is an unknown of the problem.  We will show an example (Test 4) in Section 4.  
\end{rem}
 
\subsection{Asymptotic solution of the problem}

Aim of this subsection is to study the asymptotic behavior in time of the solution of problem (\ref{model1}).
Using the result of Proposition \ref{equiv}, we will deduce that it evolves towards the unique solution $\overline{u}$  of the corresponding stationary (i.e. elliptic) obstacle problem
\begin{equation}
	\label{ell_obst}
	\overline{u}\in H^1_0(\Omega),\quad \overline{u}\geq u^c,\quad -\Delta\overline{u}\geq f,\quad (\overline{u}-u^c)(\Delta \overline u + f)=0\ \hbox{ a.e. in }\Omega.
\end{equation}

\begin{rem}\label{uniqueness}
We remark that the stationary problem corresponding to {\em(\ref{model1})}, that is
\begin{eqnarray}
\label{EQ}
	\left\{
	\arraycolsep=1.4pt\def\arraystretch{1.5}
		\begin{array}{ll}	
			H(\tilde{u}-u^c)(\Delta\tilde{u}+f)=0 &\hbox{ in } \Omega\\
			\tilde{u}=0 &\hbox{ on } \partial\Omega
		\end{array}\right.
\end{eqnarray}
is not well posed, since uniqueness fails. For example, in one dimension, if $f=0$, any function  in $H^2(\Omega)$ which coincides with the obstacle in a subset of $\Omega$ and reaches zero on the boundary in a  linear way outside of that is clearly a solution of $H(\tilde{u}-u^c)\Delta\tilde{u}=0$ and a potential asymptotic solution for problem (\ref{model1}).
%We could say ({\it \textcolor{red}{da controllare}}) that the solution of (\ref{ell_obst}) coincides with the maximal among these solutions.

Let us consider, for example, $\Omega = \left[-1,1\right]$, $f=0$ and the obstacle given by
\begin{eqnarray}
\label{u_c-esempio}
	u^c(x) = \frac 1 2 -\left(2x^2-\frac 1 2\right)^2\ ;
\end{eqnarray}
then the following functions $u_1(x)$  and $u_2(x)$ are both solutions of  problem {\em (\ref{EQ})}, but only the first is  solution of  problem {\em (\ref{ell_obst})}:

\begin{eqnarray}
\label{soluz-esempio}
	u_1(x) = \left\{
	\arraycolsep=1.4pt\def\arraystretch{1.5}
		\begin{array}{ll}
			a\left(1+x\right)	\ &	\text{\em in } -1 \leq x < -b  \\			
			u^c(x)		&	\text{\em in } -b < x \leq -0.5	\\
			0.5				&	\text{\em in } -0.5 < x < 0.5		\\		
			u^c(x)		&	\text{\em in }\ 0.5 \leq x < b			\\
			a\left(1-x\right)	&	\text{\em in } \ b \leq x \leq 1	
		\end{array}
\right.,\ \
	u_2(x) = \left\{
	\arraycolsep=1.4pt\def\arraystretch{1.5}
		\begin{array}{ll}
			a\left(1+x\right)	\ &	\text{\em in } -1 \leq x \leq -b	\\
			u^c(x)			&	\text{\em in } -b < x < b		\\
			a\left(1-x\right)	&	\text{\em in } \ b \leq x \leq 1 
		\end{array}
\right.
\end{eqnarray}
 with $a,b\in\R$ such that
$$a(1-b)=u^c(b)=u^c(-b),\quad a=(u^c)'(-b)=-(u^c)'(b).$$
 (see {\em Fig. \ref{F0}}). Note that $u^c(x), u_1(x), u_2(x)\in C^1\left(\overline{\Omega}\right)\cap H^2\left(\Omega\right)$. 
\begin{figure}[!h]
	\centering      
	\includegraphics[width=5.5cm]{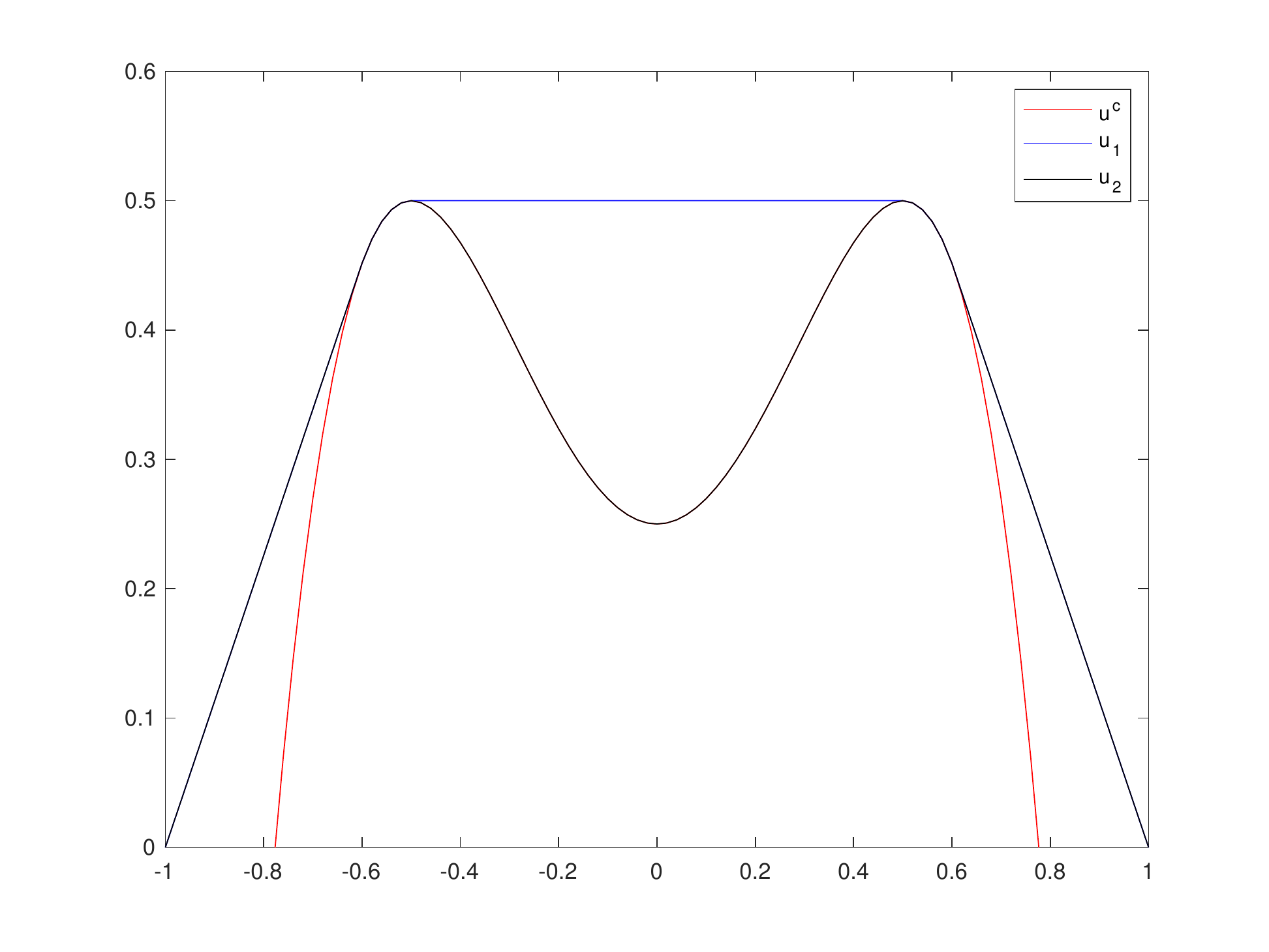}
	\caption{}
	\label{F0}
\end{figure}
If we use $u_2$ as initial datum for both problems {\em(\ref{model1})} and {\em(\ref{ell_obst})} (then with {\bf H}$_{1}$ violated), we could see that $w(t)$ would evolve in time towards $u_1$, while $u(t)$ would remain equal to $u_2$ for any time. 
\end{rem}

\begin{theor}\label{main} 
Assume conditions   \eqref{obst}, {\bf H}$_{1}$ and {\bf H}$_{2}.$ Let $u(t)$ be the global (in time) solution to the degenerate parabolic  problem  {\em(\ref{model1})}  and let $\overline{u}$  be the unique solution of the obstacle problem  {\em(\ref{ell_obst})}. Then, $u(t)$ converges to $\overline{u}$ strongly in $H^1_0(\Omega)$
for $t\to\infty,$ and there is a constant $C > 0$ such that, for every $t \geq 1$,
\begin{equation}
\label{stima1}
	\left\Vert u(t) - \overline{u}\right\Vert_{H^1\left(\Omega\right)} \leq e^{-Ct} \ .
\end{equation}
\end{theor}

\vspace{5mm}
\noindent{\bf Proof.}   
Let us make a change of variables: if we set $v(x,t)=u(x,t)-u^c(x)$, then it is easy to see that $v$ solves, for all $t\in\left(0,T\right)$, the problem
\begin{eqnarray}
\label{prob_tre}
	\left\{
	\arraycolsep=1.4pt\def\arraystretch{1.5}
		\begin{array}{ll}
			v_t - H\left(v\right)\left(\Delta v + F\right) = 0	&	\ \text{a.e. in } \Omega		\\
			v\left(0\right) =v^0						&	\text{ in } \Omega		\\
			v = g									&	\text{ on } \partial\Omega
		\end{array}
	\right.
\end{eqnarray}
with $F = \Delta u^c+f$,  $v^0=u^0 - u^c> 0$ and $g = -u^c\geq 0$ on $\partial\Omega$. For Proposition \ref{equiv} we know that $v$ also solves  the variational inequality 
\begin{eqnarray}
\label{parab_obst_v}
	\left\{
	\arraycolsep=1.4pt\def\arraystretch{1.5}
		\begin{array}{ll}
			\displaystyle v(t)\in \mathcal{K}_g, \quad \int_\Omega (v_t - \Delta v - F)( \varphi - v)\,dx \geq 0	&	\quad\forall\varphi \in\mathcal{K},\forall t\in\left(0,T\right)	\\
			v\left(0\right) = v^0	\in 	\mathcal{K}_g												
		\end{array}
	\right.
\end{eqnarray}
where
\begin{equation}
\label{convex2}
	\mathcal{K}_g= \left\{\varphi \in H^1\left(\Omega\right), \varphi - g \in H^1_0\left(\Omega\right), \varphi \geq 0 \text{ in } \Omega\right\}.
\end{equation}
The corresponding elliptic obstacle problem is then 
\begin{equation}
	\label{ellv_obst}
	\displaystyle \bar{v}\in \mathcal{K}_g, \quad  \int_\Omega(- \Delta \bar{v} - F)(\varphi - \bar{v})\, dx \geq 0 \quad	\text{for all } \varphi \in\mathcal{K}_g\ ,
\end{equation}
and $\bar{v}$ minimizes in $\mathcal{K}_g$ the following functional $$\mathcal{F}(\varphi) = \frac{1}{2}\int_\Omega\left\vert\nabla \varphi\right\vert^2 \  \mathrm{d}x- \int_\Omega F\varphi\  \mathrm{d}x.$$
We note that $-\nabla\mathcal{F}(\varphi) = \Delta \varphi + F.$

\medskip The proof follows the idea of Theorem  1.7  of \cite{CSV}. The main difference is that in \cite{CSV} it is considered $F=-1$, while here we will consider a generic datum $F$ with $F\leq 0$ 
{by  {\bf H}$_{2}.$}

First of all we prove the following  constrained Lojasiewicz inequality for the obstacle problem  (see Proposition 4.1 of \cite{CSV}) : there is a dimensional constant $C_d>0$
such that
\begin{equation}
	\label{stima2}
	\left(\mathcal{F}\left(v\right) - \mathcal{F}\left(\overline{v}\right)\right)_{+}^{\frac{1}{2}} \leq C_d\left\Vert\nabla\mathcal{F}\left(v\right)\right\Vert_{\mathcal{K}_g}
\end{equation}
for every $v\in  H^2(\Omega)\cap\mathcal{K}_g$ where
\begin{equation}
	\label{stima3}
	\left\Vert\nabla\mathcal{F}\left(v\right)\right\Vert_{\mathcal{K}_g} := \sup\left\{0,\sup_{\varphi\in\mathcal{K}_g\setminus\left\{v\right\}}\frac{-\displaystyle\int_\Omega\left(\varphi - v\right)\nabla\mathcal{F}\left(v\right)\  \mathrm{d}x}{\left\Vert \varphi - v\right\Vert_{L^2}}\right\}.
\end{equation}

We point out that the unique solution $\overline{v}$ of the obstacle problem  (\ref{ellv_obst}) solves
 $-\Delta\overline{v} = F\chi_{\left\{\overline{v} > 0\right\}}$ in $\Omega$ and $\overline{v} =g$  on $\partial\Omega.$ Then, taking $\varphi=\overline{v}$,
 
\begin{eqnarray}
	\label{stima4}
	\arraycolsep=1.4pt\def\arraystretch{3.0}
		\begin{array}{lcl}
		\left\Vert\nabla\mathcal{F}\left(v\right)\right\Vert_\mathcal{K}	&	\geq	&	-\dfrac{\displaystyle\int_\Omega\left(\overline{v} - v\right)\nabla\mathcal{F}\left(v\right)\  \mathrm{d}x}{\left\Vert\overline{v} - v\right\Vert_{L^2}} = \\
												&	=	&	-\dfrac{1}{\left\Vert\overline{v} - v\right\Vert}_{L^2}\displaystyle\int_\Omega\left(v - \overline{v}\right)\left(\Delta v + F\right)\  \mathrm{d}x = 				\\
												&	=	&	-\dfrac{1}{\left\Vert\overline{v} - v\right\Vert}_{L^2}\displaystyle\int_\Omega\left(v - \overline{v}\right)\left(\Delta v - \Delta\overline{v}\right) \  \mathrm{d}x+	\\
												&		&	-\dfrac{1}{\left\Vert\overline{v} - v\right\Vert}_{L^2}\displaystyle\int_{\Omega\cap\left\{\overline{v} = 0\right\}}\left(v - \overline{v}\right) F \  \mathrm{d}x=	\\
												&	=	&	\dfrac{1}{\left\Vert\overline{v} - v\right\Vert}_{L^2}\left(\displaystyle\frac{1}{2}\int_\Omega\left\vert\nabla\left(v - \overline{v}\right)\right\vert^2 \  \mathrm{d}x- \displaystyle\int_{\Omega\cap\left\{\overline{v} = 0\right\}}F \left(v - \overline{v}\right)\  \mathrm{d}x\right)		
	\end{array}
\end{eqnarray}

As
\begin{eqnarray}
	\label{stima5}
	\arraycolsep=1.4pt\def\arraystretch{3.0}
		\begin{array}{lcl}
		\displaystyle\int_\Omega\left\vert\nabla\left(v - \overline{v}\right)\right\vert^2	\  \mathrm{d}x&	=	&	\displaystyle\int_\Omega\left\vert\nabla v\right\vert^2 \  \mathrm{d}x+ \int_\Omega\left\vert\nabla\overline{v}\right\vert^2  \mathrm{d}x - 2\int_\Omega\nabla v\nabla\overline{v} \  \mathrm{d}x=												\\
																	&	=	&	\displaystyle\int_\Omega\left\vert\nabla v\right\vert^2 \  \mathrm{d}x+ 2\int_\Omega\left\vert\nabla\overline{v}\right\vert^2 \  \mathrm{d}x- 2\int_\Omega\nabla v\nabla\overline{v}\  \mathrm{d}x- \int_\Omega\left\vert\nabla\overline{v}\right\vert^2 \  \mathrm{d}x=	\\
																	&	=	&	\displaystyle\int_\Omega\left\vert\nabla v\right\vert^2 \  \mathrm{d}x+ 2\int_\Omega\nabla\overline{v}\cdot\left(\nabla\overline{v} \  - \nabla v\right)\mathrm{d}x - \int_\Omega\left\vert\nabla\overline{v}\right\vert^2 \  \mathrm{d}x=							\\
																	&	=	&	\displaystyle\int_\Omega\left\vert\nabla v\right\vert^2\  \mathrm{d}x - \int_\Omega\left\vert\nabla\overline{v}\right\vert^2 \  \mathrm{d}x- 2\int_\Omega\Delta\overline{v}\left(\overline{v} - v\right)\  \mathrm{d}x =											\\
																	&	=	&	\displaystyle\int_\Omega\left\vert\nabla v\right\vert^2 \  \mathrm{d}x- \int_\Omega\left\vert\nabla\overline{v}\right\vert^2\  \mathrm{d}x - 2\int_{\Omega\cap\left\{\overline{v} > 0\right\}} F\left(v - \overline{v}\right)\  \mathrm{d}x
	\end{array}
\end{eqnarray}
estimate (\ref{stima4}) becomes
\begin{eqnarray}
	\label{stima6}
	\arraycolsep=1.4pt\def\arraystretch{3.0}
		\begin{array}{lcl}
		\left\Vert\nabla\mathcal{F}\left(v\right)\right\Vert_{\mathcal{K}_g}	&	\geq	&	\dfrac{1}{\left\Vert\overline{v} - v\right\Vert}_{L^2}\left(\displaystyle\frac{1}{2}\int_\Omega\left\vert\nabla v\right\vert^2\mathrm{d}x- \frac{1}{2}\int_\Omega\left\vert\nabla\overline{v}\right\vert^2\mathrm{d}x+ \right.	\\
															&		&	- \left.\displaystyle\int_{\Omega\cap\left\{\overline{v} > 0\right\}} F\left(v - \overline{v}\right)\mathrm{d}x - \displaystyle\int_{\Omega\cap\left\{\overline{v}= 0\right\}}F \left(v - \overline{v}\right)\mathrm{d}x\right) = 	\\
															&	=	&	\dfrac{1}{\left\Vert\overline{v} - v\right\Vert}_{L^2}\left(\mathcal{F}\left(v\right) - \mathcal{F}\left(\overline{v}\right)\right)\ .
	\end{array}
\end{eqnarray}

By Poincar\'{e} inequality and    by   {\bf H}$_{2}$,  we obtain
\begin{eqnarray}
	\label{stima7}
	\arraycolsep=1.4pt\def\arraystretch{3.0}
		\begin{array}{lcl}
		\left\Vert v - \overline{v}\right\Vert^2_{L^2\left(\Omega\right)}	&	\leq	&	C_p\left\Vert\nabla\left(v - \overline{v}\right)\right\Vert^2_{L^2\left(\Omega\right)} =																																								\\
														&	=	&	C_p\left(\displaystyle\int_\Omega\left\vert\nabla v\right\vert^2  \  \mathrm{d}x- 2\int_{\Omega\cap\left\{\overline{v} > 0\right\}}vF\  \mathrm{d}x - \int_\Omega\left\vert\nabla\overline{v}\right\vert^2 \  \mathrm{d}x+ 2\int_{\Omega\cap\left\{\overline{v} > 0\right\}}\overline{v}F\ \mathrm{d}x\right)	\\
														&	\leq	&	2C_p\left(\mathcal{F}\left(v\right) - \mathcal{F}\left(\overline{v}\right)\right)
	\end{array}
\end{eqnarray}
and so
\begin{equation}
	\label{stima8}
	\left\Vert\nabla\mathcal{F}(v)\right\Vert_{\mathcal{K}_g} \geq \frac{1}{\left\Vert v - \overline{v}\right\Vert_{L^2}}\left(\mathcal{F}\left(v\right) - \mathcal{F}\left(\overline{v}\right)\right) \geq \frac{1}{\sqrt{2C_p}}\frac{\mathcal{F}\left(v\right) - \mathcal{F}\left(\overline{v}\right)}{\left(\mathcal{F}\left(v\right) - \mathcal{F}\left(\overline{v}\right)\right)^{\frac{1}{2}}}\ ,
\end{equation}
that is the following constrained Lojasiewicz inequality  holds 

\begin{equation}
	\label{stima9}
	\left(\mathcal{F}\left(v\right) - \mathcal{F}\left(\overline{v}\right)\right)^{\frac{1}{2}}      \leq \sqrt{2C_p}\left\Vert\mathcal{\nabla F}(v)\right\Vert_{\mathcal{K}_g}.
\end{equation}

Then, by using Proposition 2.10 in \cite{CSV}, we conclude the proof, since

\begin{equation}
	\left\Vert u - \overline{u}\right\Vert_{H^1\left(\Omega\right)} =\left\Vert v - \overline{v}\right\Vert_{H^1\left(\Omega\right)} \leq e^{-Ct} \ . 
\end{equation} \qed

\vspace{5mm}

\begin{rem}
Let us consider the following quantities: $$ M(t)=\int_\Omega (u(t)-u^c) \ dx, \quad I(t)= \int_\Omega H(u(t)-u^c)(\Delta u (t)+f) \ dx.$$ The first one measures  the global distance in time of the solution from the obstacle. Theorem \ref{main} and (\ref{EQ}) imply that 
$$ M(t)\to \overline{M}= \int_\Omega (\bar{u}-u^c) \ dx, \qquad I(t) \to 0 \ .$$
If we integrate the equation in (\ref{model1}) we get
\begin{equation}
	\frac{d}{dt} M(t) = \int_\Omega u_t(t) \ dx = \int_\Omega H(u(t)-u^c)(\Delta u (t)+f) \ dx = I(t) .
\end{equation}
 When $I(t)$ does not change its sign in time, then the convergence of $M(t)$ is monotone. To look at the time profiles of $M(t)$ and $I(t)$ is  interesting, since it gives some informations on the global evolution of the solution (for example, it reveals the contact times with the obstacle). We will look at their corresponding discrete quantities in the numerical simulations of the last section. 
\end{rem}

%%%%%%%%%%%%%%%%%%%%%%%%%%%%%%%%%%%%

\section{Numerical approximation}

We start for simplicity with the one dimensional setting, and $\Omega=(-1,1)$. Then the problem to solve becomes:
\begin{eqnarray}
	\label{Bm}
		\left\{
		\arraycolsep=1.4pt\def\arraystretch{1.5}
		\begin{array}{ll}
			u_t -H(u-u^c)(u_{xx}+f)=0		&\text{ in } \Omega\times  (0,T) 		\\ 
			u(0)=u^0	 				&\text{ in } \Omega 				\\ 
			u=0 						&\text{ on } \partial\Omega\times  (0,T), 
		\end{array}\right.
\end{eqnarray}
where $T$ is a sufficiently large time, $u^0 > u^c$ and $H$ denotes the  Heaviside function defined in (\ref{heavi}), or eventually a regular approximation  of it in a small right neighborhood of the origin, for example the $C^1$ function $\eta_n$ given by:

\begin{equation}
	\label{eta1}
	\eta_n(r)=\left\{
	\arraycolsep=1.4pt\def\arraystretch{1.5}
		\begin{array}{ll} 
			1				& \text{ if } r>\frac1 n 		\\ 
			-2{n^3r^3}+3n^2r^2	& \text{ if } 0\le r \le \frac1 n 	\\   
			0 				& \text{ if } r<0 
		\end{array} 
		\right. .
\end{equation}

According to the fixed integer parameter $n\in\N$, we see that $\eta_n(r)=H(r)$ for $r\ge \frac 1 n$ and $r\le 0$, and that  $$ \| H-\eta_n\|_{L^1(-1,1)}=\int_{-1}^1 |H(r)-\eta_n(r)|\ dr = \int_0^{1/n} |H(r)-\eta_n(r)|\ dr \to 0 \quad as\ n\to\infty\ .$$

Note that, as happens for $H$, $\eta_n(0)=0$, so that even in this case the diffusion coefficient vanishes at the contact points between the solution of  problem (\ref{model1}) (with $H$ replaced by $\eta_n$) and the obstacle. But now also all the values of $\eta_n$ in $(0,1]$ characterize supercritical states of the solution close to the obstacle itself.

On $\Omega$ we define for a given $N\in\N$ a uniform  grid $G_h$ of  size $ h=2/N$. Then $G_h$ will have $N-1$ internal nodes $x_j=-1+jh$ ($j=1,..,N-1$) over a total number of $(N+1)$. 
 
Concerning  time discretization, we adopted a uniform time step $ \Delta t=T/M$, for a given $M\in\N$, so that  the  solution is computed at any time $t^k=k\Delta t$ ($k=1,..,M$):  by $u^k_{j}$ we denote the discrete solution at time  $t^{k}$ in a node $x_{j}$ of $G_h$. The initial data will be given by
\begin{equation}
	\label{id}
	u_{j}^{0}=u^0(x_{j}), \ u_{j}^{c}=u^c(x_{j}),\text{ with }\ u_{j}^{0}> u_{j}^{c} \text{\  for any }j=1,..,N-1,\   u_{0}^{0}=u^0_N=0.
\end{equation}

We are interested in the numerical solution of  (\ref{Bm}) on the grid ${G}_{h}$;  to avoid stability problems without heavy restrictions on the parabolic step ratio $\gamma=\Delta t/h^2$ we adopted a  semi implicit finite difference scheme: 

\bigskip($S$): $ for\ any\  k=0,1,...,M-1\ solve \ for \ any \ internal \ node\ x_j :\ \hskip10cm$
\begin{eqnarray}
	\label{scheme}	
		\left\{
		\begin{array}{l}
		\begin{array}{l}
		u_{j}^{k+1}=u_{j}^{k}+\Delta t\ z_{j}^{k}\ (\delta_h u^{k+1}_{j}+f_j):= 
		u_{j}^{k}+\gamma\ z_{j}^{k}\ \left(u^{k+1}_{j-1}-2u^{k+1}_{j} +u^{k+1}_{j+1}\right)+\Delta t\ z_{j}^{k}\ f_j\ , \\  \\ \text{if } u_{j}^{k+1}< u_{j}^{c} \Rightarrow u_{j}^{k+1}=u_{j}^{c}\ ,\\ \\
		{u_{0}^{k+1}=u_{N}^{k+1}=0 \quad  (\text {boundary values})}	;	\end{array}
		\end{array}
		\right.
\end{eqnarray}
with $\delta_h$ we have indicated  the usual 3-point second order finite difference operator over $G_h$; we will talk of scheme ($S_H$) when $z^k_j=H^k_j:=H\left(u_{j}^{k}-u_{j}^{c}\right)$ (that is when we use the sharp Heaviside values), of scheme ($S_\eta$) when $z^k_j=\eta^k_j:=\eta\left(u_{j}^{k}-u_{j}^{c}\right)$ (that is when we use its approximated values given by (\ref{eta1})). Then at any time iteration $k$ one has to solve the following linear system:  $$ B^ku^{k+1}:=(I +\gamma z^k *A) u^{k+1} = u^k +\Delta t z^{k}F\ ,$$ where $u^k$, $z^k$ and $F=(f_j)$ are column vectors of dimension $(N-1)$, $A$ is the  tridiagonal $(N-1)\times(N-1)$ matrix associated to the discrete Laplacian in  one dimension, and by $v*M$ we  mean the vector matrix product in which    each line $j$ of $M$ is multiplied for the $ j-$th component of $v$.

\begin{rem} 
Let us explain the second line of scheme (\ref{scheme}).  We  proved in Section 2 that the solution of (\ref{model1}) always remains over the obstacle. In the discrete settings with scheme (\ref{scheme}) anyway, the impact with the obstacle happens at a certain time iteration, with a thrust which depends on the  parameter $\gamma$ and which can cause the overcoming of the obstacle  before the Heaviside term can stop the diffusion. So it is necessary to force the discrete solution to coincide with the obstacle where it has gone over. 
When $\gamma$ is large, anyway, the solution can overcome the obstacle in many adjacent nodes at a single instant time, yielding an overestimation of the contact set which the subsequent iterations are no more able to correct. That is why, even if the scheme has no stability constraints, a reduced value of $\gamma$ (hence of $\Delta t$) should be necessary  in order to evolve towards the correct stationary solution, with a consequent grow of  computational costs. 

To face such a problem we have experimented some variants of our approach. The first one consists in the use of the approximate Heaviside function $\eta_n$ of {\em(\ref{eta1})} to determine the diffusion coefficient: when the solution gets closed to the obstacle, it has the effect to reduce progressively the thrust and even to prevent the overcome of the obstacle (if a suitable value of $n$ is chosen). 

Another idea is to use a variable discretization time step, reducing it only when it is necessary. We tested two ways for that. The first one measures the impact thrust of each $\Delta t$ in terms of the number of nodes involved in the contact at a single iteration time, halving it  until this number remains large but resetting it at the initial value when the contact with the obstacle becomes sufficiently stable. It works well, but this \lq\lq trial and error" process is still too expensive. The second way comes directly from the scheme. Assume for simplicity $f=0$;  if $u^k_j\geq u^c_j$ for any $j$,  in order to remain over the obstacle everywhere at the $k+1$ iteration we should have
$$u^{k+1}_j=u^k_j+\Delta t z^k_j \delta_hu^{k+1}_j\geq u^c_j,\quad \forall j\ ;$$
where there is already a contact ($z^k_j=0$) there is nothing to prove; elsewhere $z^k_j=1$ and if the solution decreases  at a node $x_j$ then necessarily $\delta_hu^{k+1}_j<0$, so that the previous inequality is equivalent to ask

\begin{equation}
	\label{step_estimate}
	\Delta t \leq D_j:=\frac{u^k_j-u^c_j}{-\delta_hu^{k+1}_j}\ ;
\end{equation}
then the estimate of the smallest positive value of $D_j$ (with $\delta_hu^{k+1}_j$ replaced by  $\delta_hu^{k}_j$) gives at any iteration a sufficiently small time step in order to reach the obstacle with the right thrust.  We have tested  all these approaches in the experiments of the next section, trying a comparison evaluation.
\end{rem}

\bigskip\noindent In order to  emphasize the convergence of the solutions to the stationary state, as discussed in the previous section, we adopted for scheme (\ref{scheme}) the following stopping criterium: 
%we left the schemes work in time until a suitable condition is satisfied, such as  for example
%$$\| u^k -u^{k-1}\|\infty < tol$$ 
%(that is when the solution is sufficiently stabilized), or  
\begin{equation}
	\label{stop}
	\max_j\left[(u^k_j-u^c_j)|\delta_h u^k_j+f_j|\right]< tol\ ,
\end{equation} 
where $tol$ indicates a prescribed small tolerance. In other words the scheme stops before the final time $T$ if  the limit problem is sufficiently solved.

For sake of comparison, we have also implemented a numerical scheme for the corresponding parabolic and elliptic obstacle problems (respectively (\ref{parab_cs}) and (\ref{ell_obst})), with the same discretization parameters, showing even at a discrete level the essential coincidence of  the solutions of the two evolutive problems (if $\gamma$ is not too large) and their convergence to the same asymptotic solution. Many algorithms can be found in the literature for the obstacle problem: among them we have choosen the ones presented in \cite{BS}, based on the iterative solutions of piecewise linear systems. The discrete version of the equation in (\ref{parab_cs}) becomes
	\begin{equation}
		\label{BS}
			(w^{k+1}-u^c)^T(w^{k+1}+\gamma A w^{k+1}-w^k-\Delta t f)=0\ .
	\end{equation}

Setting $y=w^{k+1}-u^c\ge0$, then $y$ has to solve $$y^T(y+\gamma Ay-b)=0\ ,$$ with $b=w^k-u^c+\Delta t f-\gamma Au^c$. In \cite{BS} it is proved that $y=\max(x,0)$ is a solution of the previous equation if $x$ solves

	\begin{equation}
		\label{PL}
			[I+\gamma AP(x)]x=b ,
	\end{equation}
where $P(x)$ is the diagonal matrix with $p_{jj}=H(x_j)$, and $H$ is the Heaviside (sign) function (\ref{heavi}). In order to solve the last implicit equation a quasi-Newton method is implemented which needs a certain number of linear system solutions (Picard iterations) for any discrete time step: 
$$P^0=O\ (\hbox{null\ matrix}), \quad (I+\gamma A P^n)x^{n+1}=b,\quad for \ n=0,1,... \ \hbox{ until } \ P^n=P^{n+1}\ ;$$ then $x=x^{n+1}$ is the solution of (\ref{PL});  hence $w^k=y+u^c$ solves (\ref{BS}) and evolves in time towards the solution $\bar{u}$ of the corresponding stationary obstacle problem (\ref{ell_obst}) on the grid $G_h$.

\medskip The extension of scheme (\ref{scheme}) to the two-dimensional case is straightforward, at least when $\Omega$ is a rectangular open set $(a,b)\times(c,d)$. Using equal space steps $\Delta x=\Delta y=h$, the discrete solution $u^k_{ij}$ will denote the approximated value of $u$ in the node $x_{ij}$ at time $t^k=k\Delta t$. It is then sufficient to replace the finite difference operator $\delta_h$ with the usual five-point Laplacian approximation scheme:
$$ \delta^2_h u^k_{ij}=\frac{u^k_{i+1j}+u^k_{i-1j}-4u^k_{ij}+u^k_{ij+1}+u^k_{ij-1}}{h^2} .$$
All the previous considerations remain unchanged.

%%%%%%%%%%%%%%%%%%%%%%%%%%%%%%%%%%%%

\section{Numerical tests}  
We have tested  scheme  (\ref{scheme}) with the stopping criterium (\ref{stop}), for $tol=10^{-4}$ and different initial data and obstacles in one dimension on $\Omega=(-1,1)$ and in two dimensions on square regions. Here we discuss the results of these experiments.

\begin{figure}[!ht]
	\centering      
	\includegraphics[width=4.5cm]{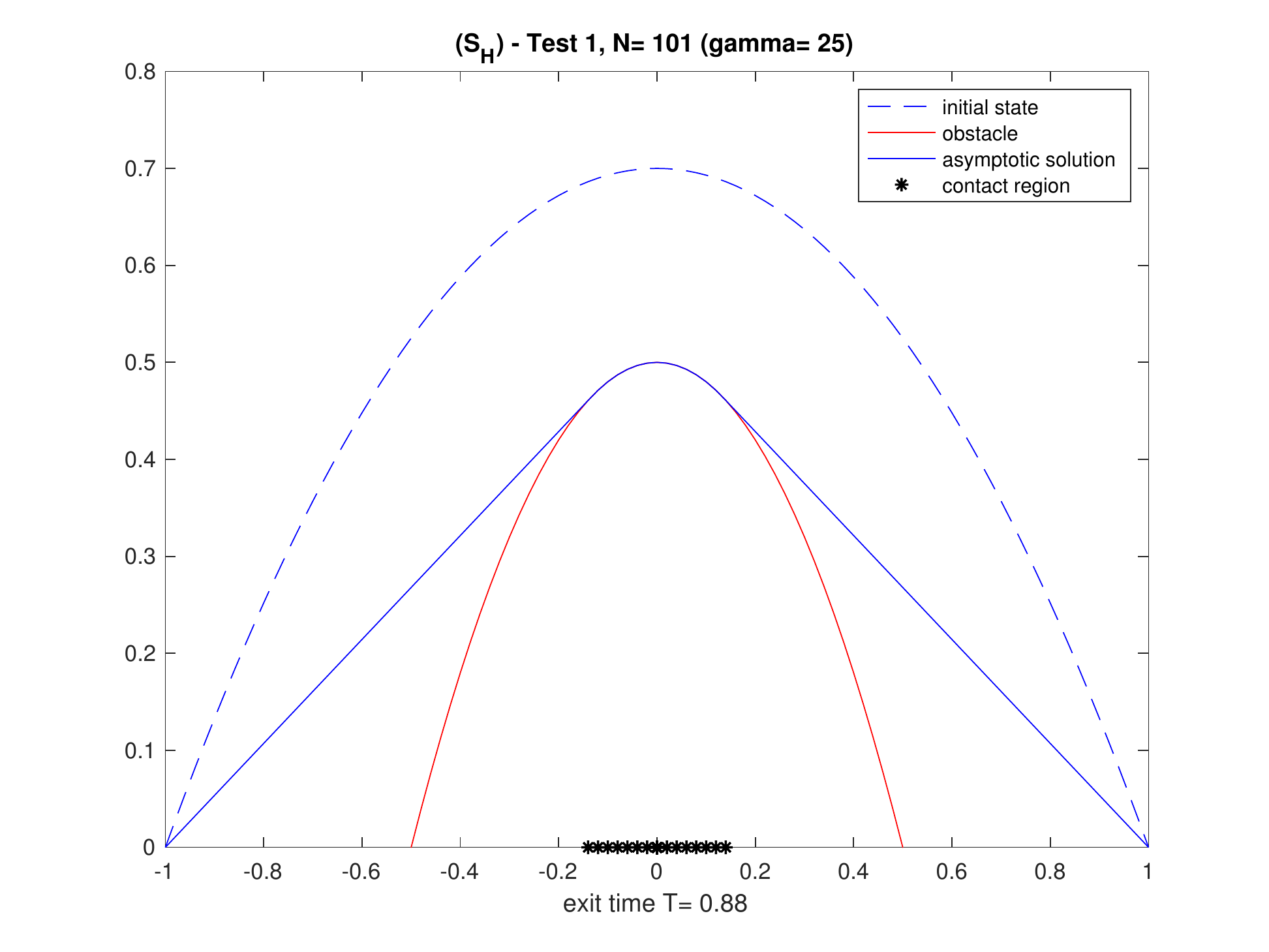}
	\includegraphics[width=4.5cm]{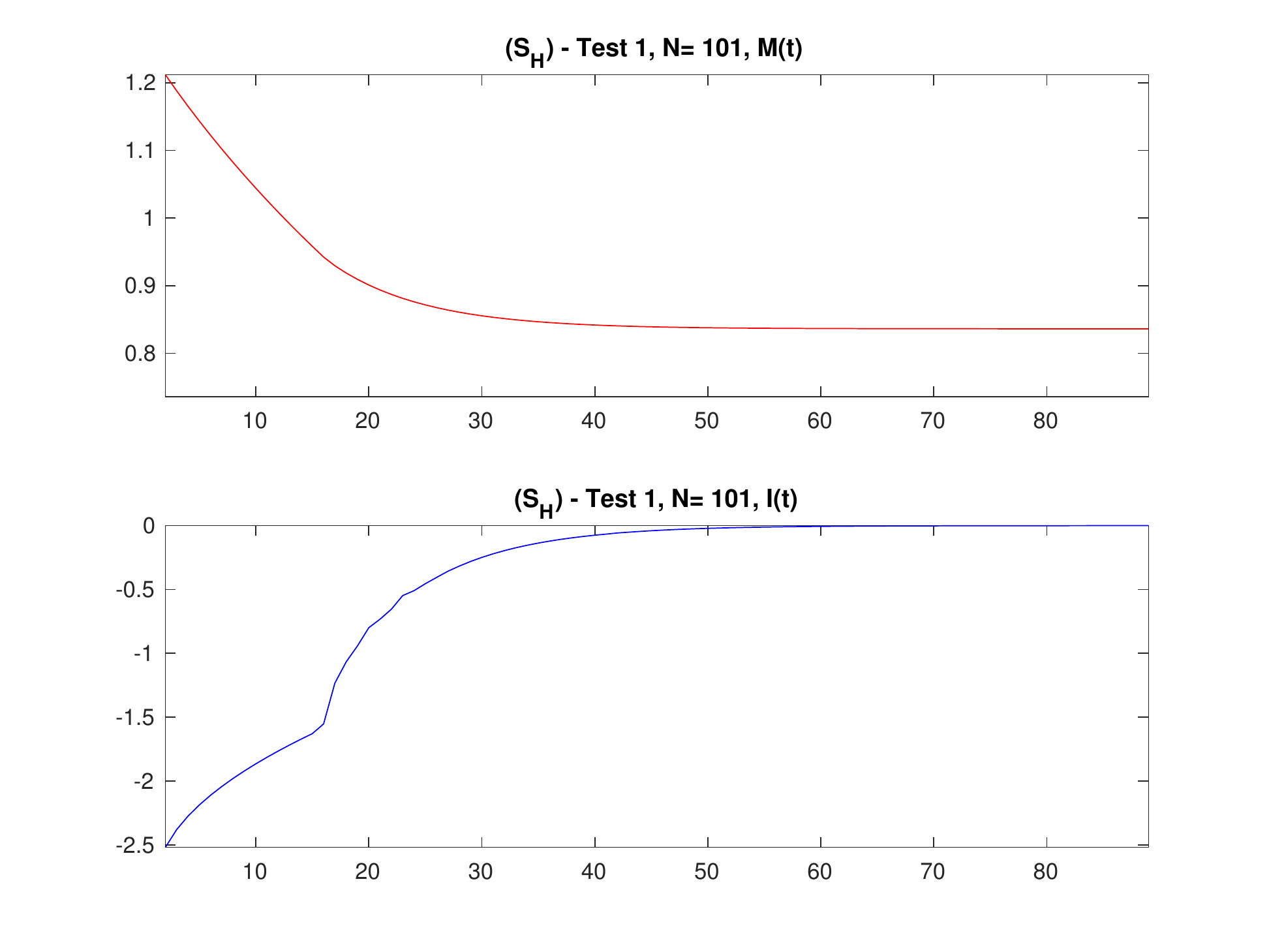}
	\includegraphics[width=4.5cm]{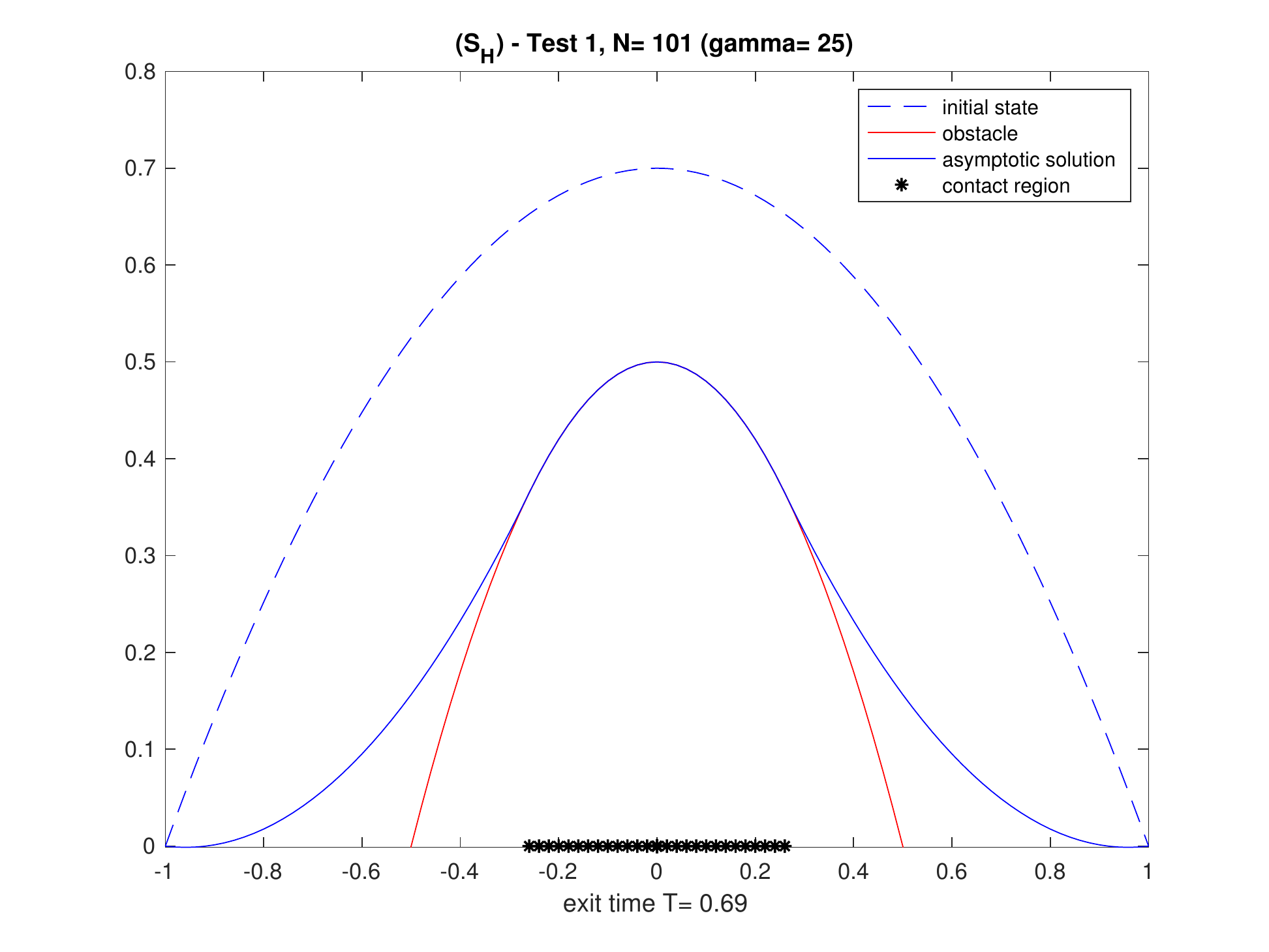}
	\caption{\footnotesize{Test 1. a) $f=0$, b) discrete $M(t)$ and $I(t)$ evolution;  c) $f=-1.5$.}	}
	\label{F1}
\end{figure}

\medskip{\bf Test 1.} $u^0=0.7-0.7x^2$, $u^c=0.5-2x^2$ (inverted parabola, with negative values at $\partial\Omega$); when $f=0$, the solution decreases in time until it touches the obstacle from the top; the two lateral branches (in the detachment region) then rapidly become linear, that is harmonic, until  nothing changes anymore (Fig.\ref{F1} a). The discrete contact region (with $N=101$ nodes) is the set $C=[-0.14,0.14]$.  In Fig.\ref{F1} b) the plots are reported of the discrete quantities corresponding to $M(t)$ and $I(t)$, which in this case are monotone in time. The impact time of the solution with the obstacle is highlighted by the change of slope in the second plot. The  addition of a constant negative source term ($f=-1.5$) correctly increases the contact set (now $C=[-0.26,0.26]$) and reduces the final stopping time (Fig.\ref{F1} c).

\begin{figure}[!ht]
	\centering      
	\includegraphics[width=4.5cm]{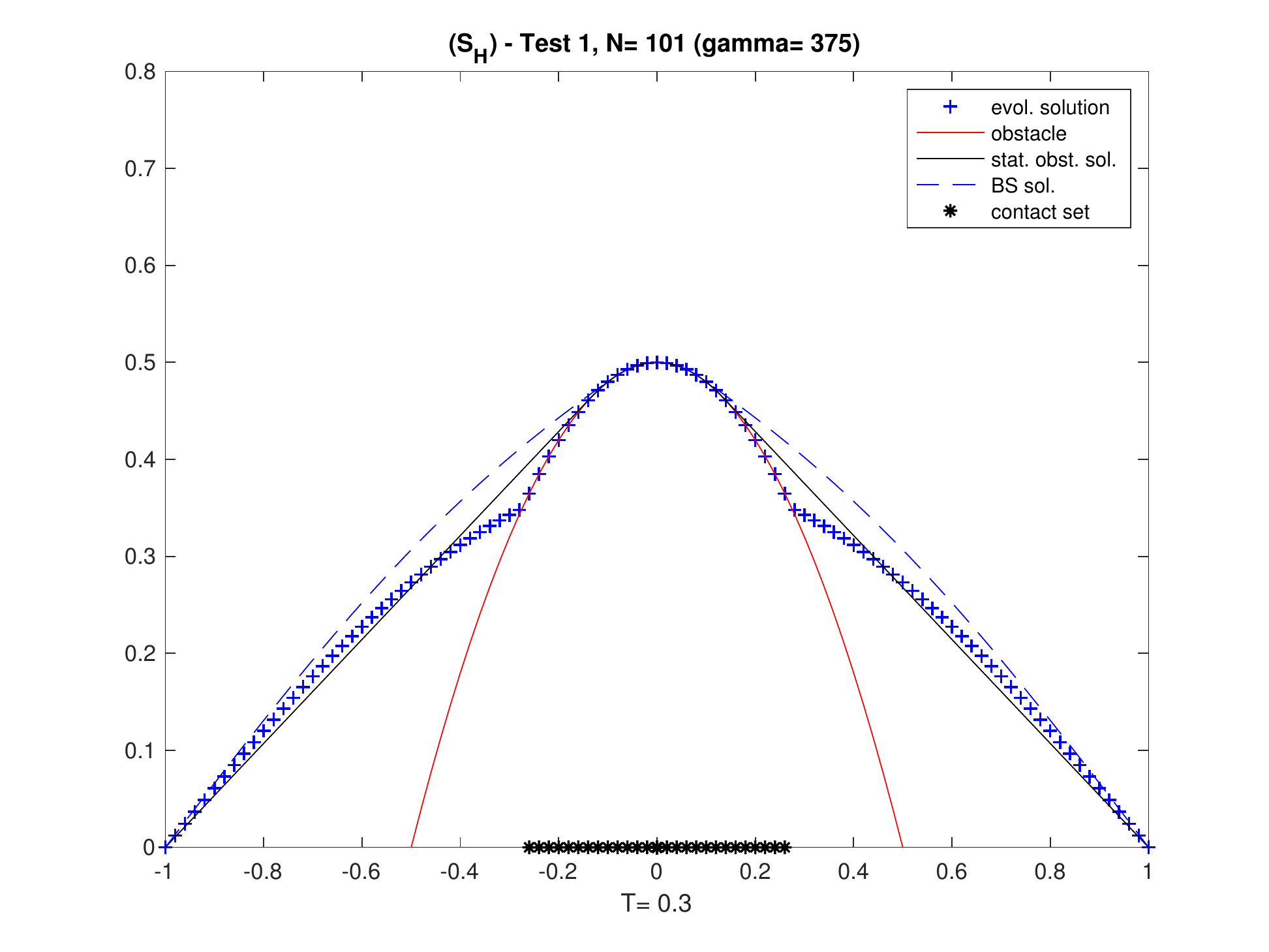}
	\includegraphics[width=4.5cm]{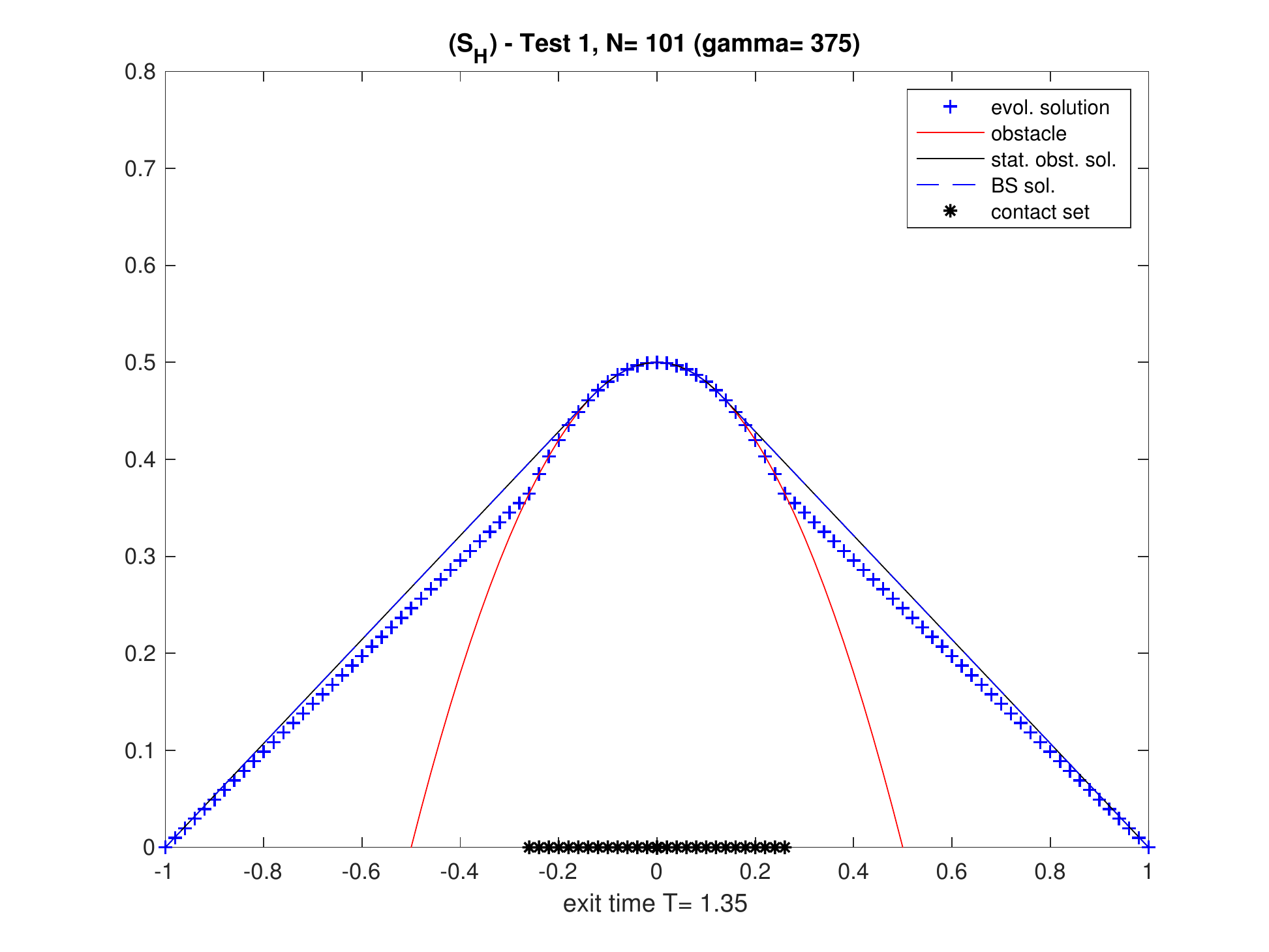}
	\caption{\footnotesize{Test 1. Overestimation of the contact set for large $\gamma$: a) first time impact; b) final uncorrect solution.}}	
	\label{F2}
\end{figure}

On this example we tried a comparison, in terms of precision and computational costs, of the different approaches introduced in the previous section. In Table 1 the first column indicates the  type of Heaviside function (H=exact, $\eta_n$=approximated), the second one if a fixed (F) or variable (V) time step approach (the one based on the time step estimate (\ref{step_estimate})) is adopted during the evolution; $T^*$ denotes the exit time reached applying  criterium (\ref{stop}), $C_{bound}$ the right extremum of the detected symmetric contact set with the obstacle (which as we said should be in this case $0.14$), $\|u-w\|_\infty$ the maximum norm of the difference in time between the discrete solutions of schemes  (\ref{scheme}) and (\ref{BS}), that is:
$$ \|u-w\|_\infty=\max_{k}\| u^k-w^k \|_{\infty} \ .$$
The table values show with a certain evidence some aspects of the different approaches:

\begin{itemize}
	\item  if $\gamma$ is too high, the contact set  can be overestimated, and the asymptotic solution is incorrect (see also Fig.\ref{F2});
	\item in order to have a good coincidence between $u$ and $w$, a low value of  $\gamma$ is necessary, that is a little $\Delta t$ and many time iterations;
	\item the use of an approximated Heaviside function with a sufficiently low parameter $n$ helps a little, since the right contact set can be found, and a better coincidence between the two solutions in time. But the evolution is slowed down in an artificial way, and the contact is less sharp;
	\item a better performance comes from the variable step approach, where, without a significative change in the exit time, the correct solution and contact set are recovered. A higher number of time iterations is needed, but much less than the one needed (with a consistent reduction of $\gamma$) in order to get  the same precision.
	\item the table also allows a cost comparison between our semi implicit approach to the obstacle problem and the implicit one of (\ref{BS}): while in the first one for any time step a single linear system has to be solved, in the second one a certain number of linear system solutions is needed.  For example, with $\gamma=75$ at the end the total number of these resolutions is of the order of 400, much more than the total time iterations of scheme (\ref{scheme}), even in its  variable time step version. Consequently, our approach to the numerical resolution of problem (\ref{model1})  can be  considered as a competitive algorithm for the approximation of the parabolic variational inequality (\ref{parab_obst}). 
\end{itemize}
\begin{table}[!h]
\begin{center}
\begin{tabular}{|c|c|c|c|c|c|c|}
\hline
Heav & time step & $\gamma$  &$T^*$  & time iter. &   $C_{bound}$ & $\|u-w\|_\infty$  \cr
\hline
H & F & 375  & 1.35 & 10  & {\bf0.26} & $6\ 10^{-2}$ \cr
\hline
H & V & 375  & 1.35 & 28  & 0.14 & $1.2\ 10^{-3}$ \cr
\hline
H & F &187.5  & 1.05 & 15  & {\bf0.2} & $3.4\ 10^{-2}$ \cr
\hline
$\eta_{20}$ & F &187.5  & 1.275 & 18  & 0.14 & $1.25\ 10^{-2}$ \cr
\hline 
H & V &187.5  & 1.12 & 34  & 0.14 & $6\ 10^{-4}$ \cr
\hline 
H & F &150  & 1.08 & 19  & 0.14 & $1.4\ 10^{-2}$ \cr
 \hline
H & F &75  & 0.96 & 33  & 0.14 & $1.4\ 10^{-2}$ \cr
 \hline
$\eta_{50}$ & F &75  & 1.56 & 53  & 0.14 & $4.1\ 10^{-3}$ \cr
 \hline
H & V &75  & 0.96 & 50  & 0.14 & $2.3\ 10^{-4}$ \cr
 \hline
H & F &37.5  & 0.9 & 61  & 0.14 & $1.8\ 10^{-4}$ \cr
\hline
H & F &18.75  & 0.86 & 116  & 0.14 & $4.4\ 10^{-4}$ \cr
\hline
H & F &9.37  & 0.84 & 226  & 0.14 & $6.6\ 10^{-4}$ \cr
\hline
\end{tabular}
\caption{\footnotesize{Test 1. Performance comparison of scheme (S) with exact (H) or approximated ($\eta_n$) Heaviside function, fixed (F) or variable (V) time step. }}
\end{center}
\label{default}
\end{table}

\begin{figure}[!ht]
	\centering      
	\includegraphics[width=3.6cm]{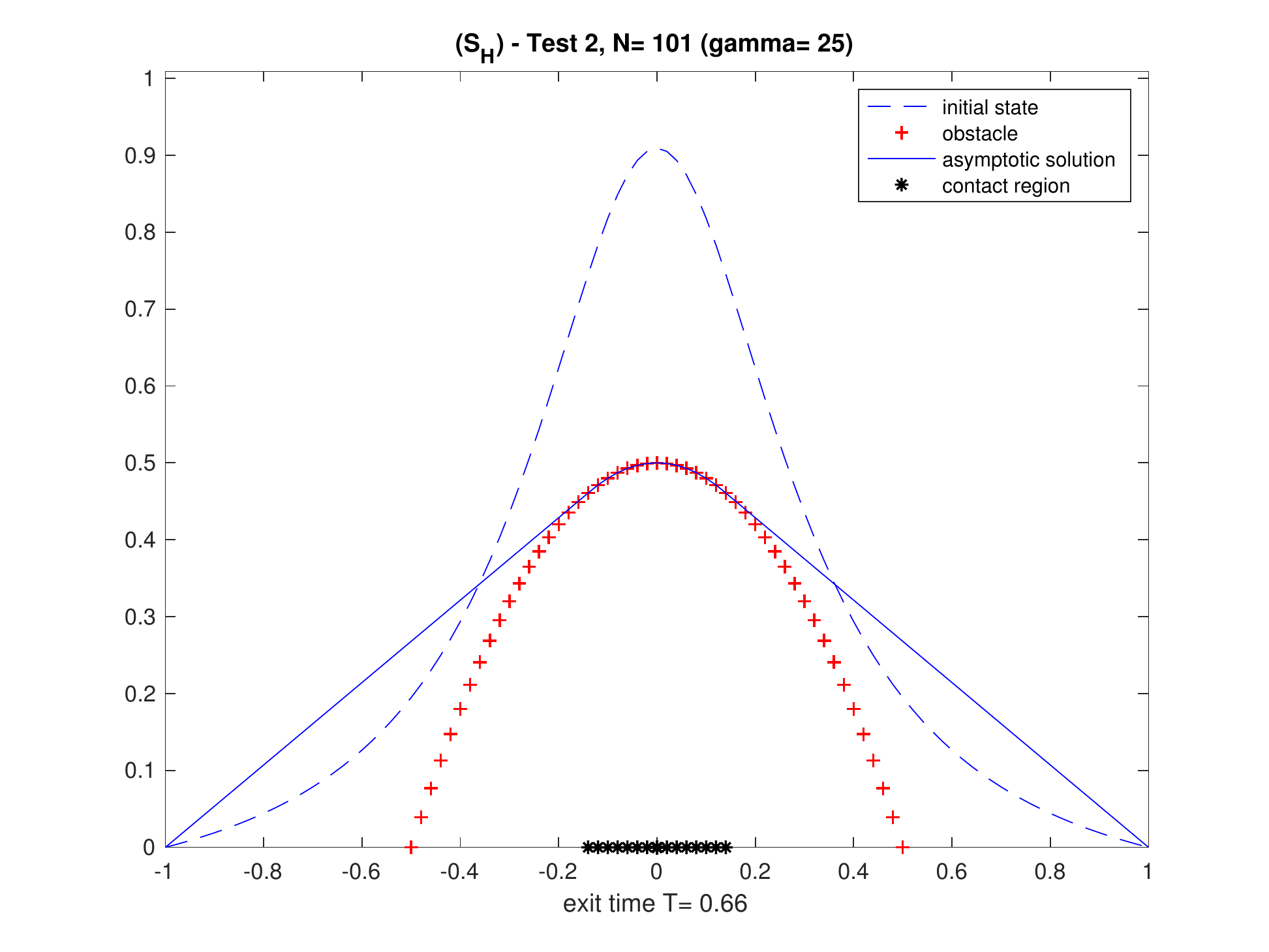}
	\includegraphics[width=3.6cm]{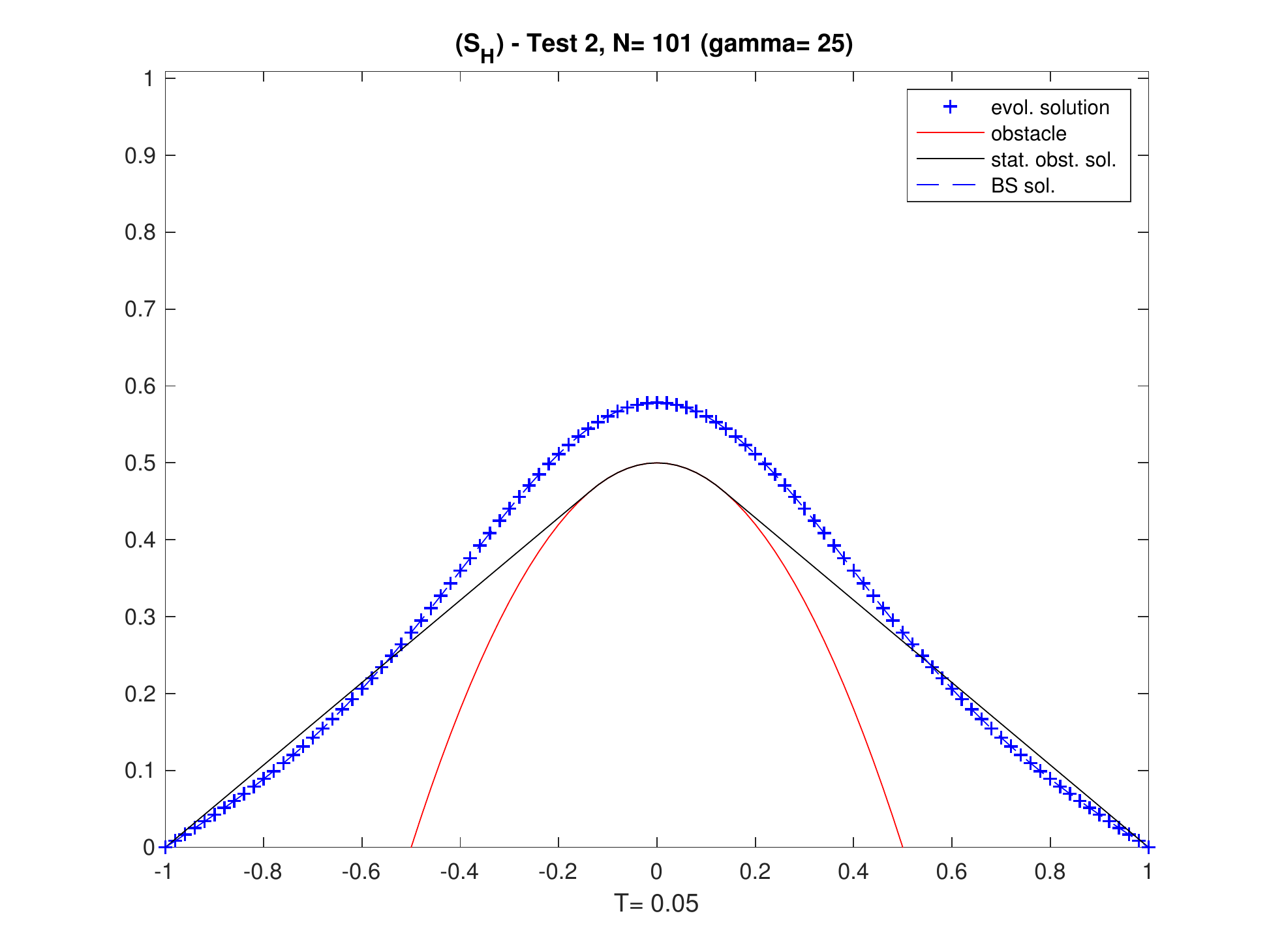}
	\includegraphics[width=3.6cm]{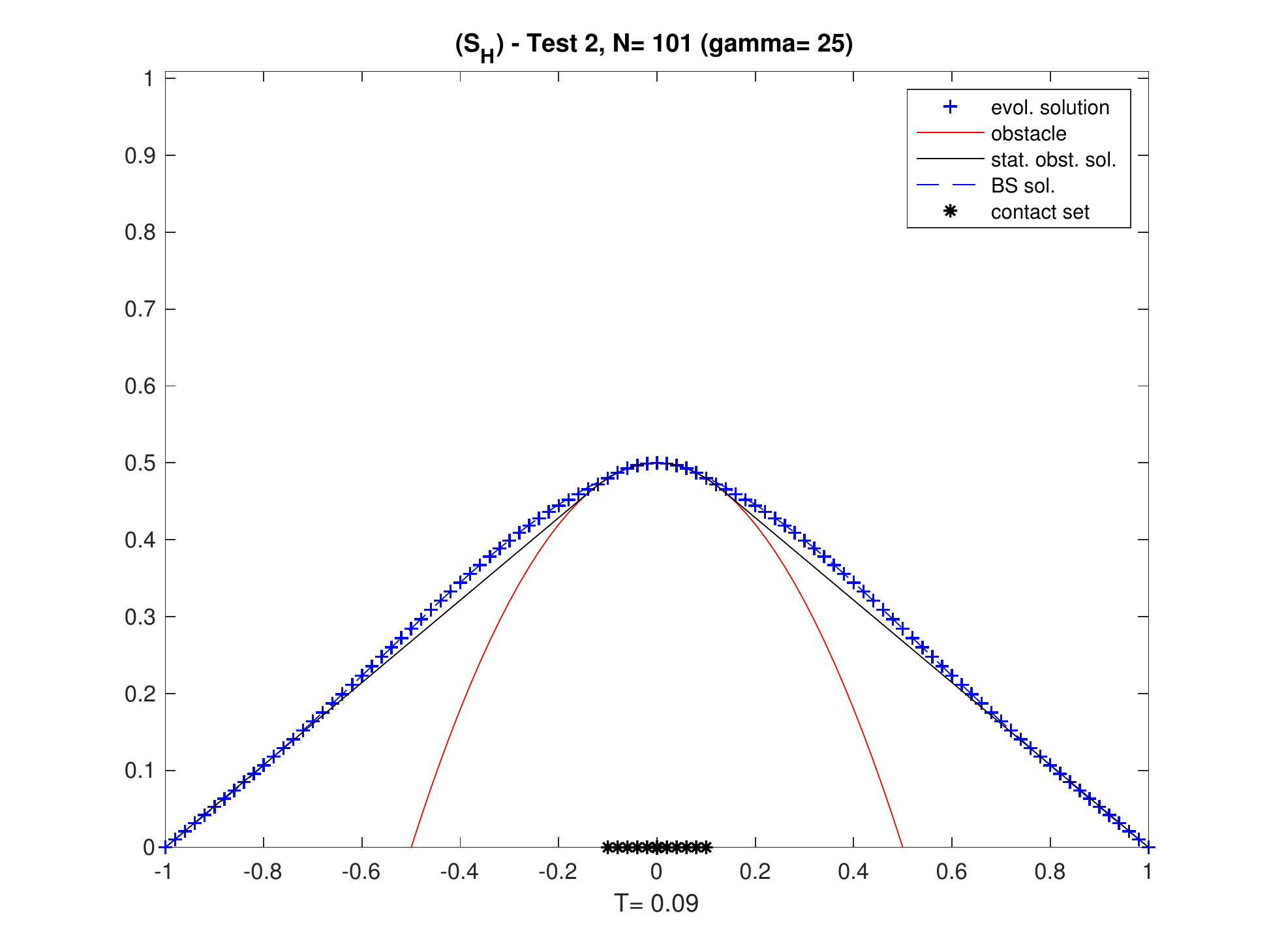}
	\includegraphics[width=3.6cm]{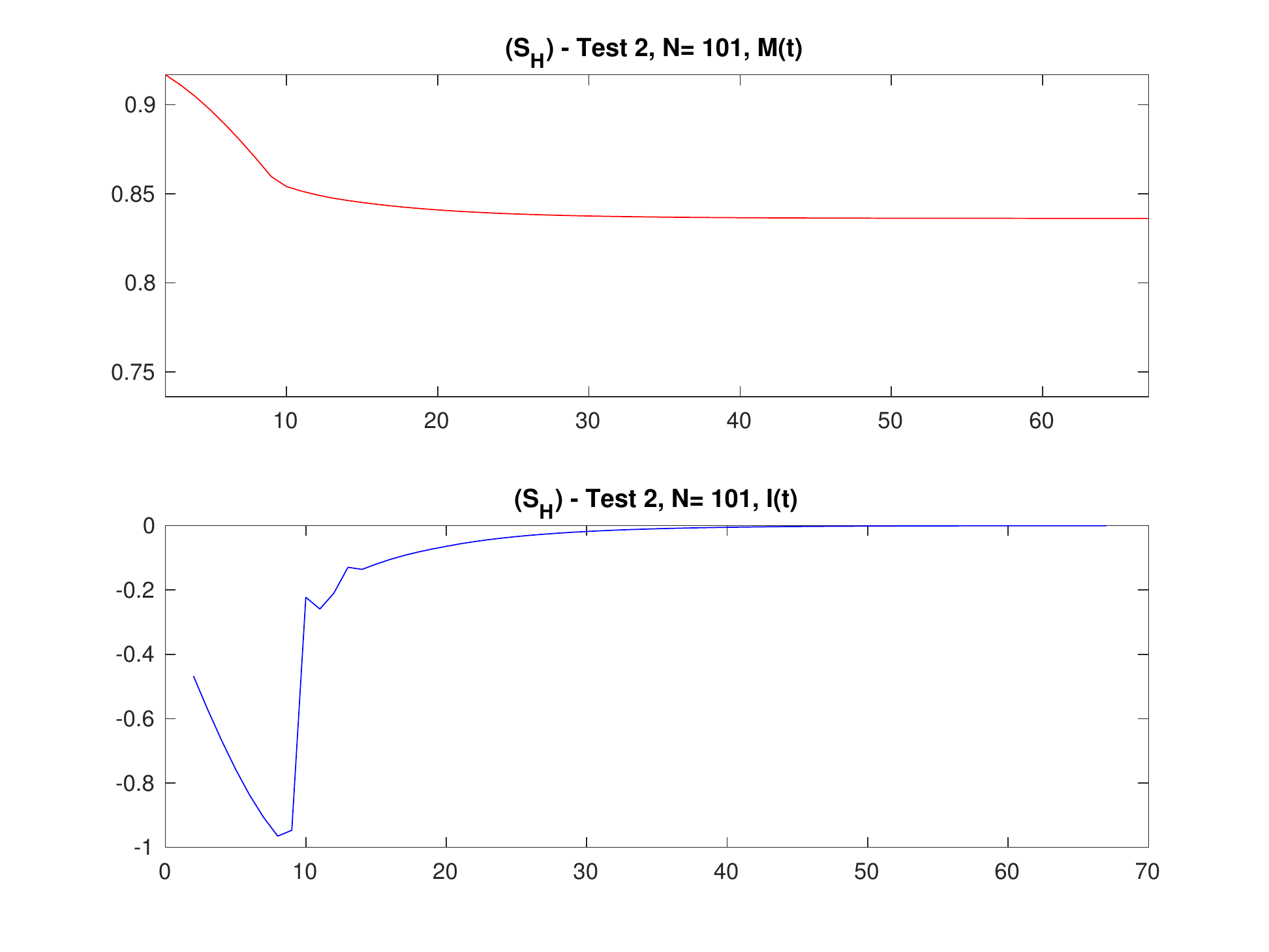}
	\caption{\footnotesize{Test 2. a) initial datum and final solution, b) t=0.05, c) t=0.09, d) discrete $M(t)$ and $I(t)$ evolution.}}	
	\label{F3}
\end{figure}

\medskip{\bf Test 2.}  $u^0=\frac 1{(1+10x^2)}-\frac 1 {11}$ (partially convex initial state), same obstacle and source term of Test 1; we get the same stationary solution of Test 1, but a different evolution (Fig.\ref{F3}). Note that now the solution initially grows in regions where it is convex and decreases where it is concave: despite of that, the total mass $M(t)$ decreases for any time. On the contrary, the quantity $I(t)$ decreases during the first part of evolution, before increasing towards zero, remaining all the time negative.

\begin{figure}[h]
	\centering      
	\includegraphics[width=3.6cm]{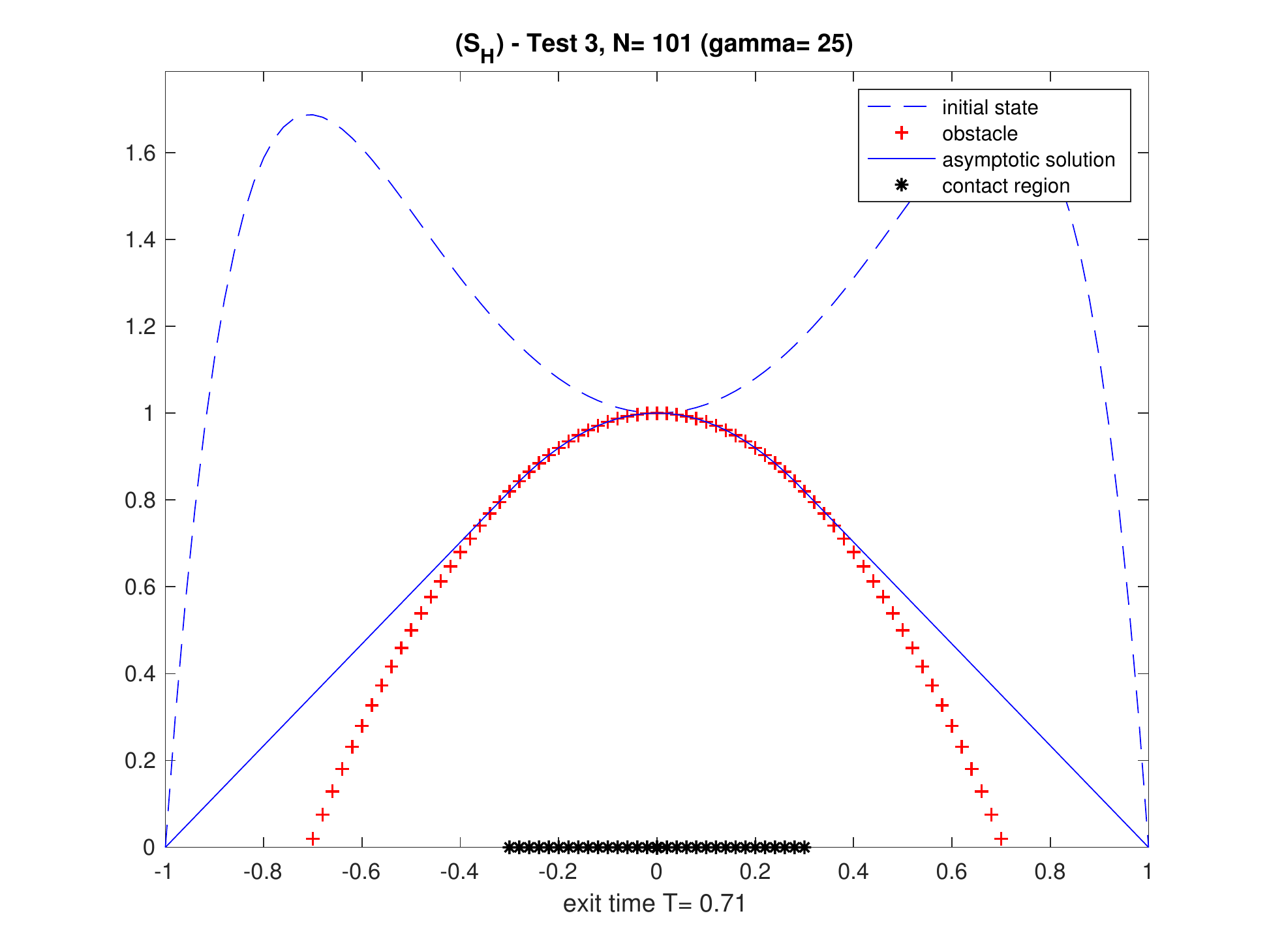}
	\includegraphics[width=3.6cm]{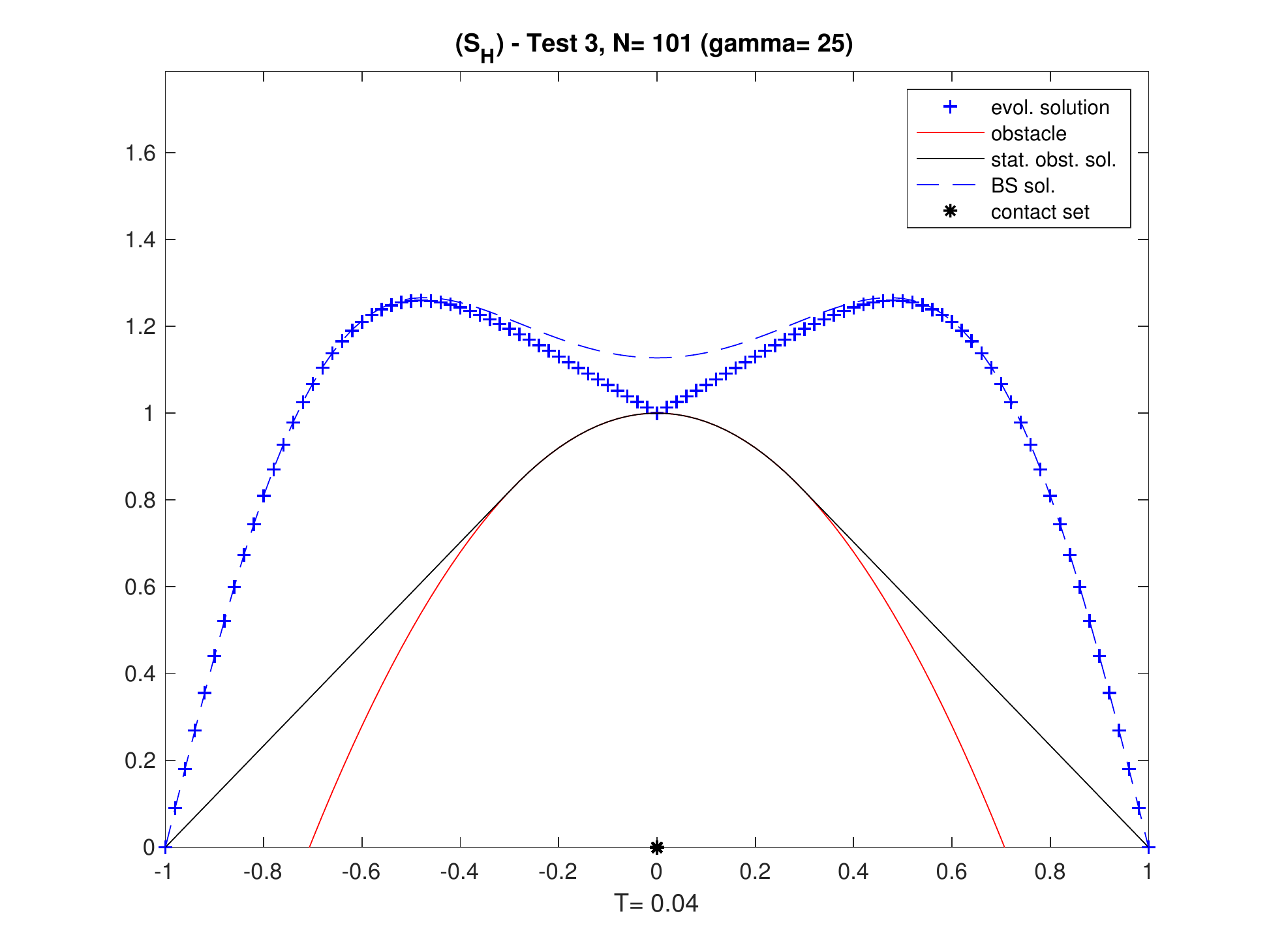}
	\includegraphics[width=3.6cm]{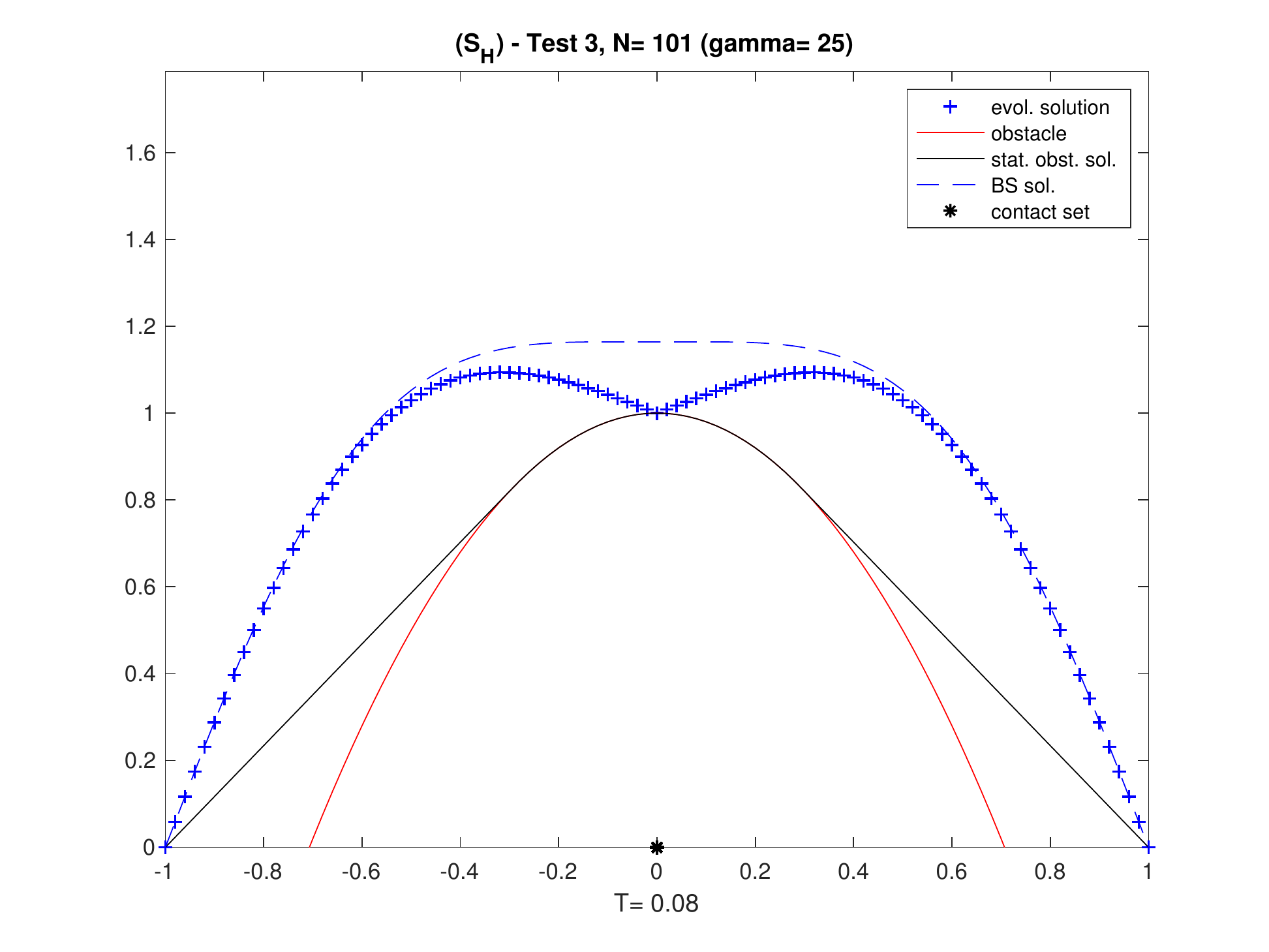}
	\includegraphics[width=3.6cm]{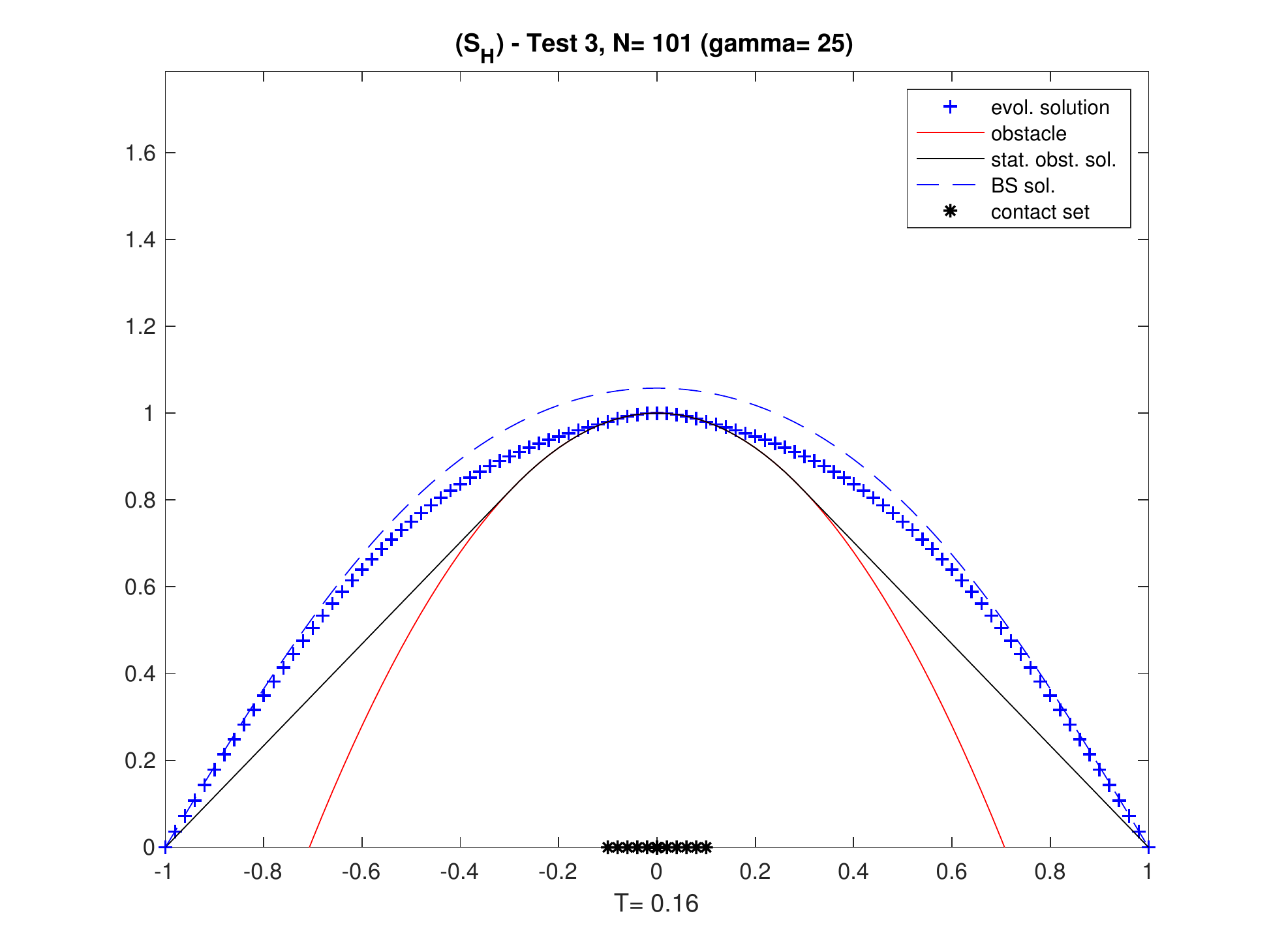}
	\caption{\footnotesize{Test 3. a) initial datum and final solution, b) t=0.04, c) t=0.08, d) t=0.16.}}
	\label{F5}
\end{figure}

\medskip{\bf Test 3.} $u^0=(1-x^2)(1+x^2)^3$, $u^c=1-2x^2$ (initial contact point with the obstacle at the origin), $f=0$; this example shows that the assumption $u^0>u^c$ is essential in order to have the same evolution of the corresponding parabolic obstacle problem (see  Remark \ref{u0}). Here the asymptotic solution is the same for the two problems, and even the final contact set is the same ($C=[-0.3,0.3]$), but the evolution is completely different: in the contact point the solution of (\ref{scheme}) (++) cannot detach anymore from the obstacle, differently to what happens to the other one (dotted), see Fig.\ref{F5}. 

\begin{figure}[!h]
	\centering      
	\includegraphics[width=4.5cm]{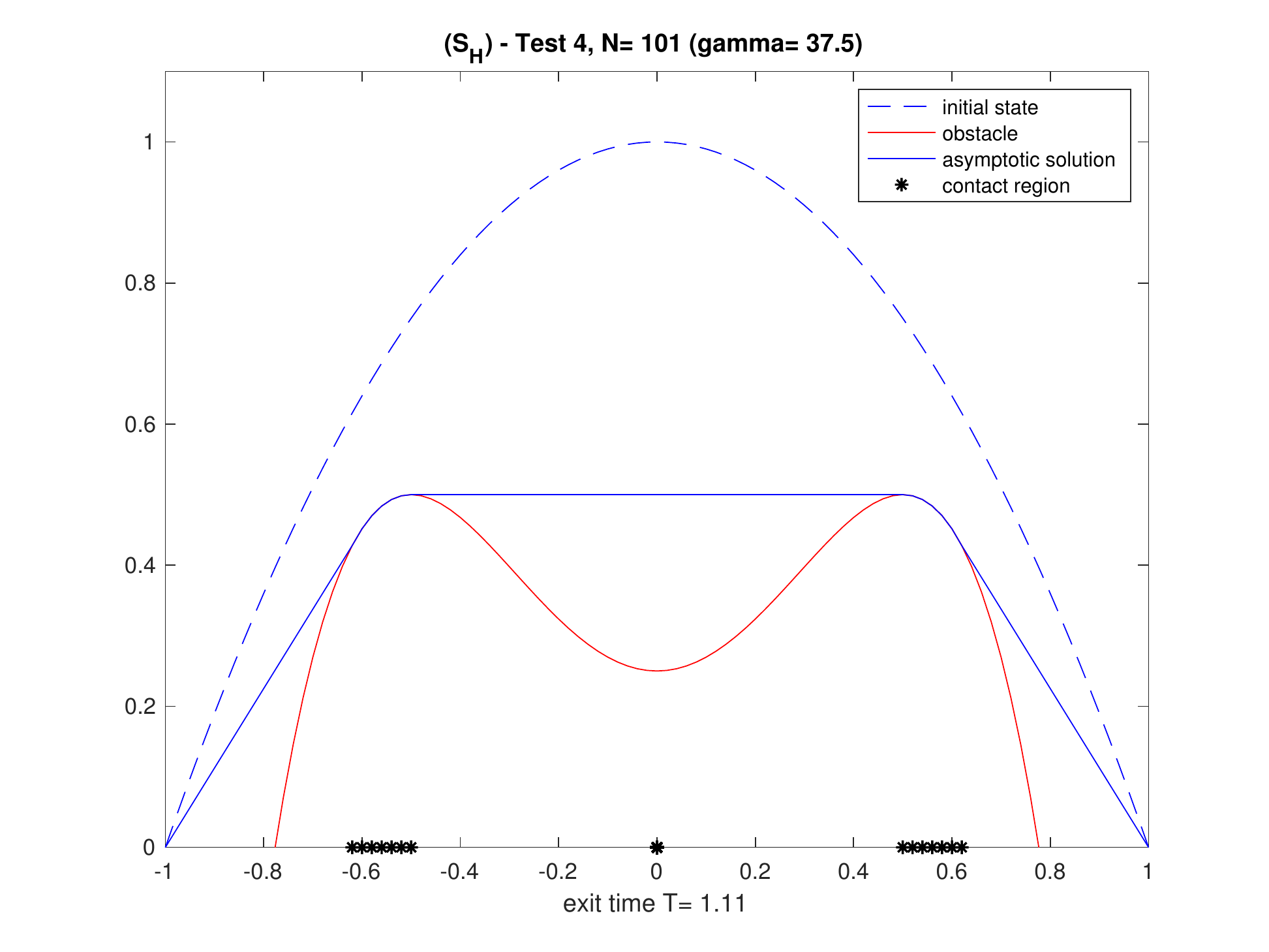}
	\includegraphics[width=4.5cm]{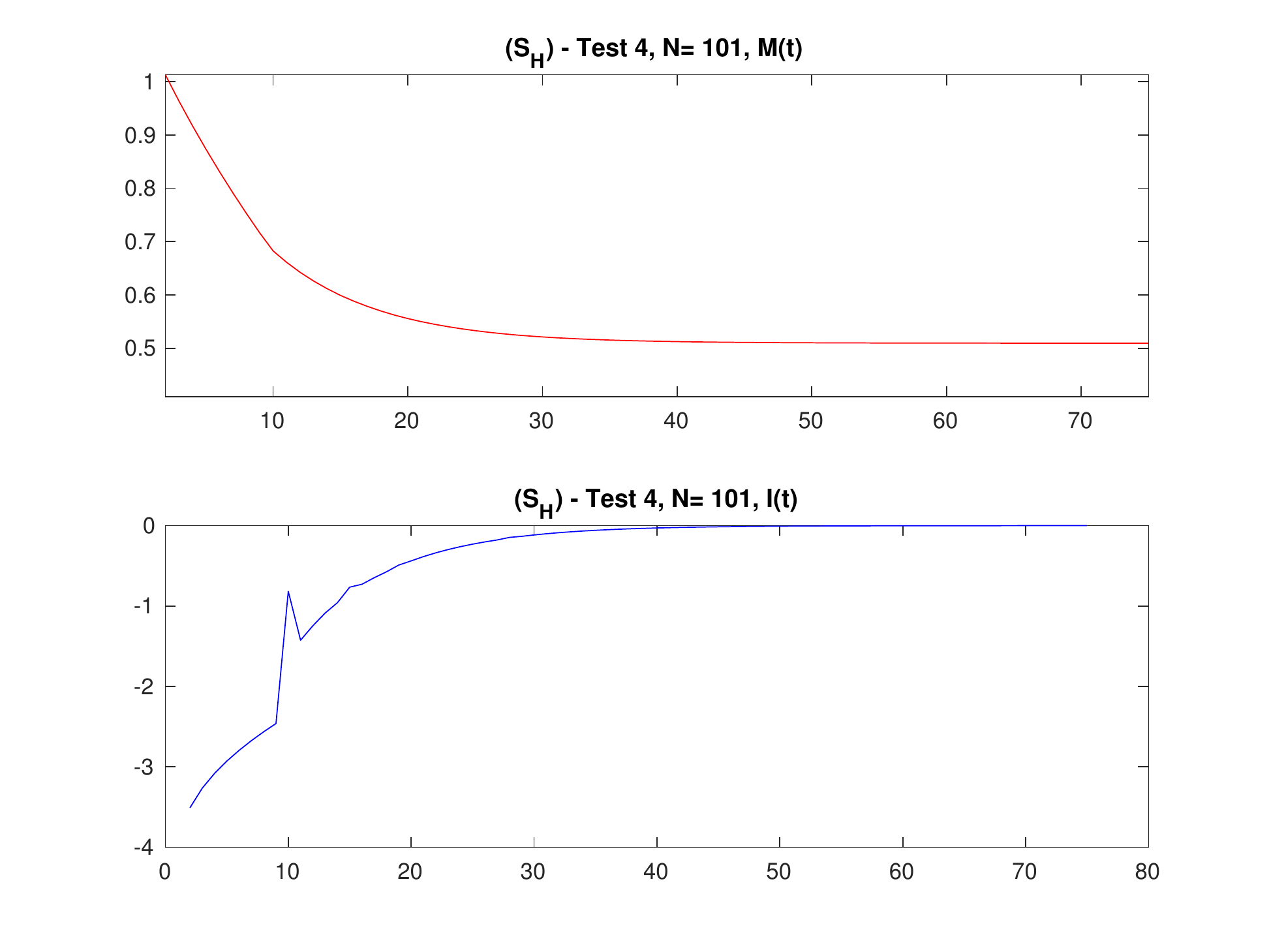}
	\includegraphics[width=4.5cm]{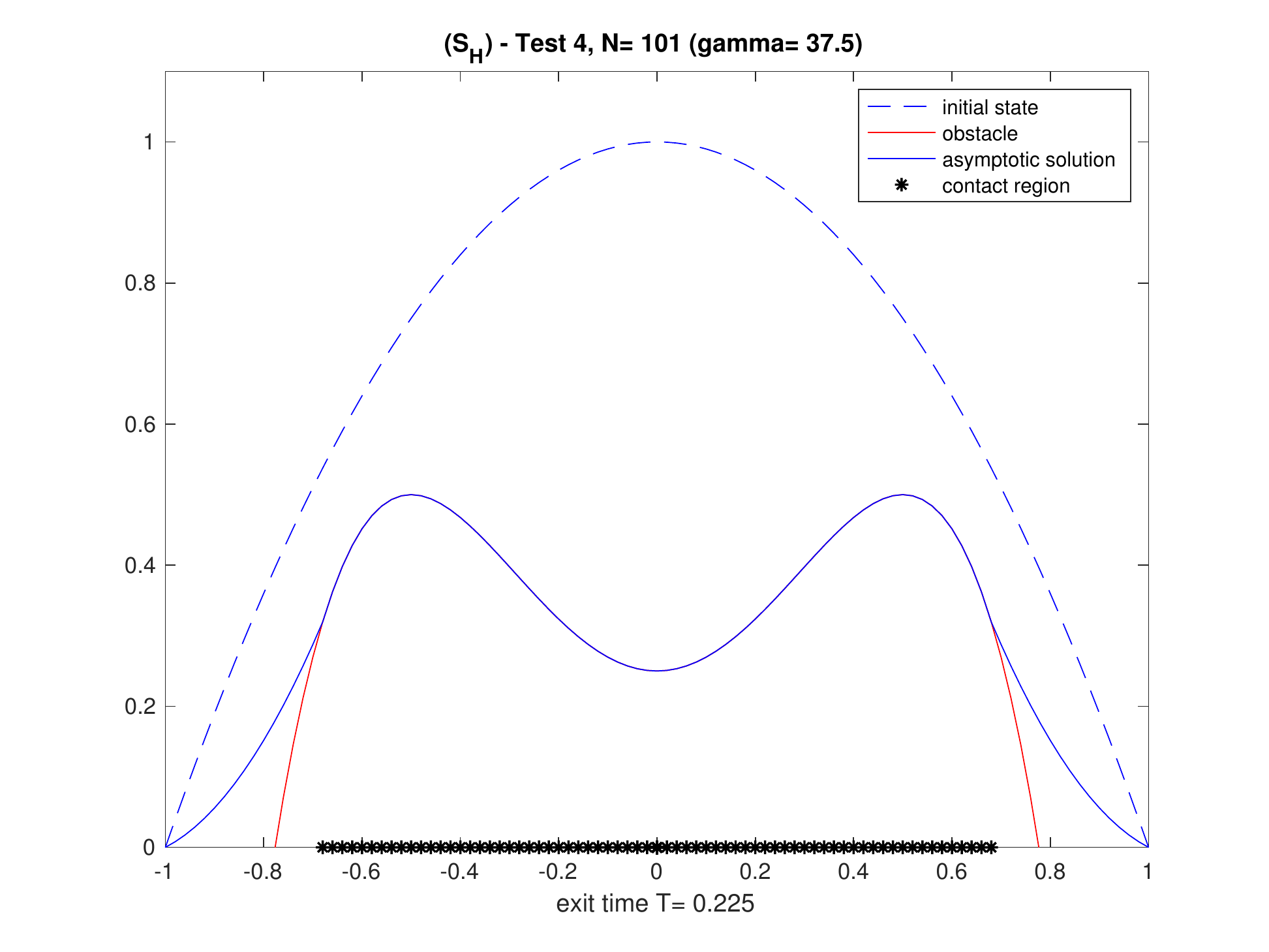}
	\includegraphics[width=4.5cm]{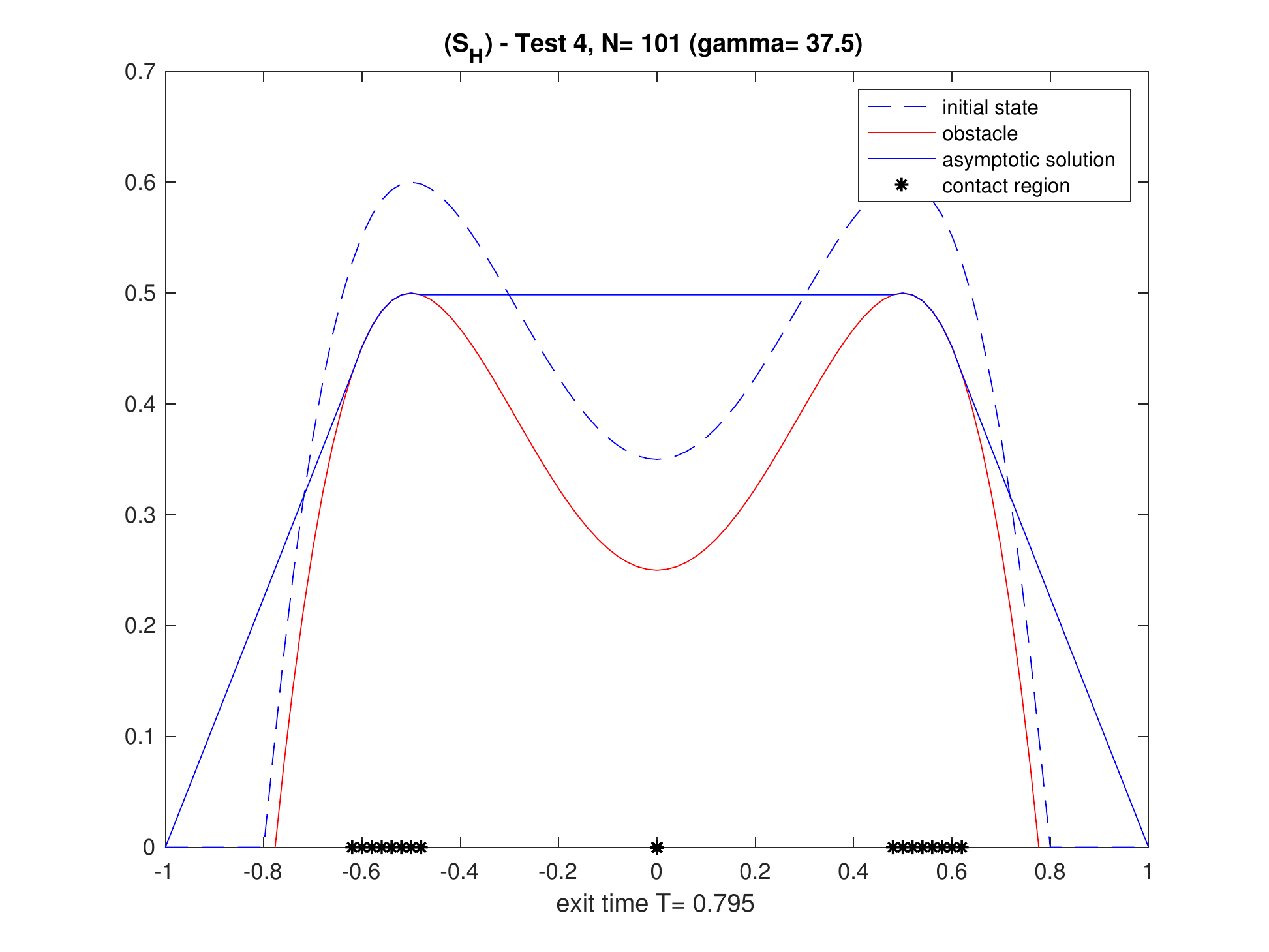}
	\includegraphics[width=4.5cm]{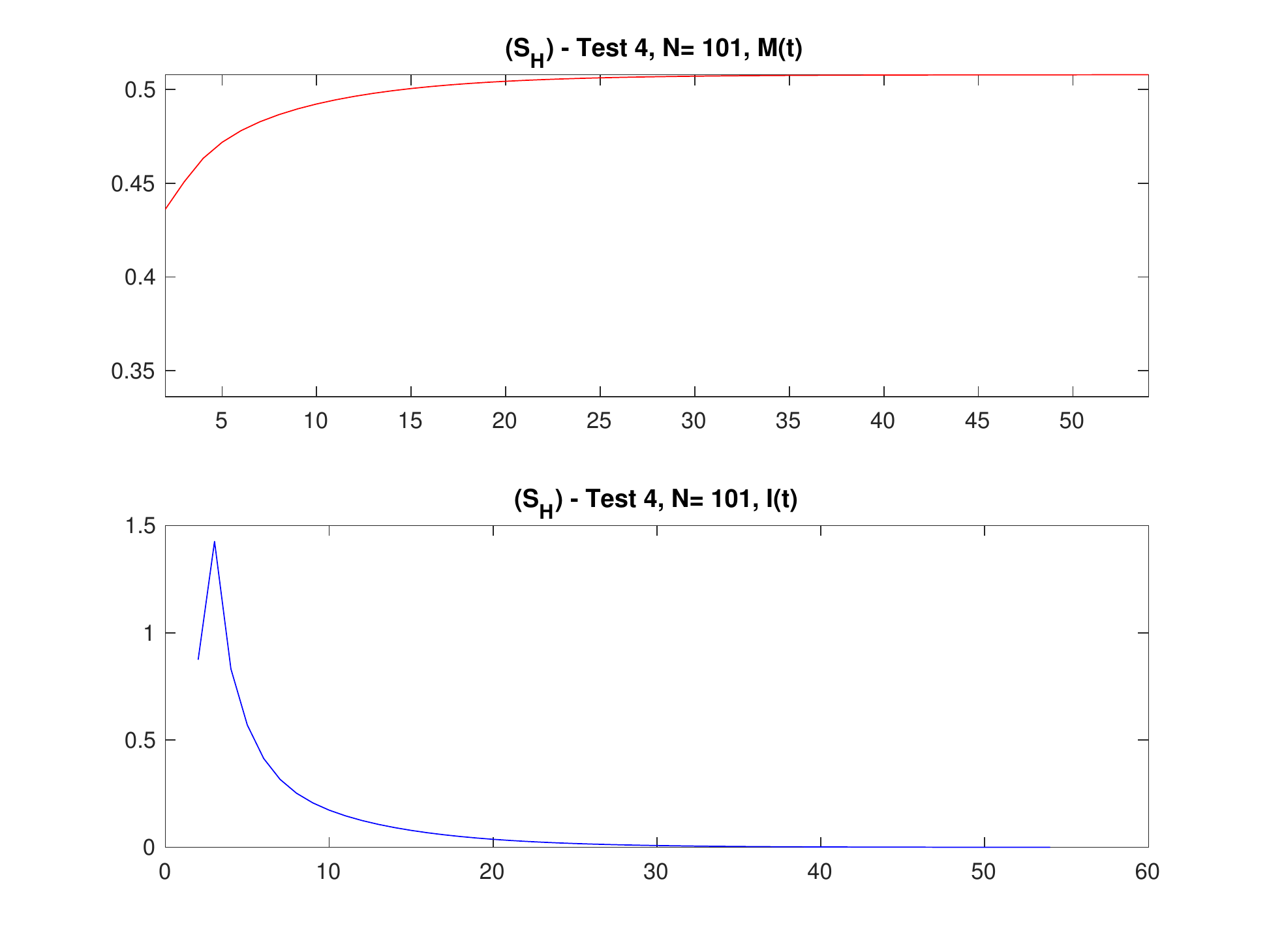}
	\caption{\footnotesize{Test 4. a) $f=0$, initial datum and final solution, b)  $M(t)$ and $I(t)$ evolution, c) $f=-4$,\\ d-e) $f=0$ but  initial datum close to the obstacle. }}	
\label{F4}
\end{figure}

\medskip{\bf Test 4.} $u^0=1-x^2$, $u^c=0.5-(2x^2-0.5)^2$ (two equal hills with a valley in the middle): it is the example of Remark \ref{uniqueness}. When $f=0$ the solution leans on the hills and remains stretched over the valley (Fig.\ref{F4} a). The final contact region is now given by $C=(-b,-0.5)\cup(0.5,b)$, with $b\simeq 0.6054$. Note that assumption {\bf H}$_{2}$  in this case is not satisfied in a small neighborhood of the origin which does not belong to the contact set. Even in this case $M(t)$ and $I(t)$ are monotone (Fig.\ref{F4} b). In order to push the solution in contact with the whole convex region of the obstacle a sufficiently negative source term has to be added: in this case $f=-4$ is necessary to make $\Delta u^c+f\leq 0$, so that  {\bf H}$_{2}$ holds in all $\Omega$, and in particular in the whole connected contact region $C=(-0.66,0.66)$  (Fig.\ref{F4} c).

 In Fig.\ref{F4} d-e we show what happens if we start with a different initial datum very close to the obstacle:
$$ u^0=\max(0,0.5-(2x^2-0.5)^2+0.1) \ ;$$
the solution converges towards the same asymptotic solution, but now essentially from below; then $I(t)$ tends to zero from positive values and $M(t)$ is monotone increasing.

\medskip  In the next two examples we considered less regular obstacles, not differentiable or even discontinuous. The experiments show that  model (\ref{model1}) still works also in these cases and that the scheme (\ref{scheme}) behaves correctly.

\begin{figure}[!h]
	\centering      
	\includegraphics[width=4.5cm]{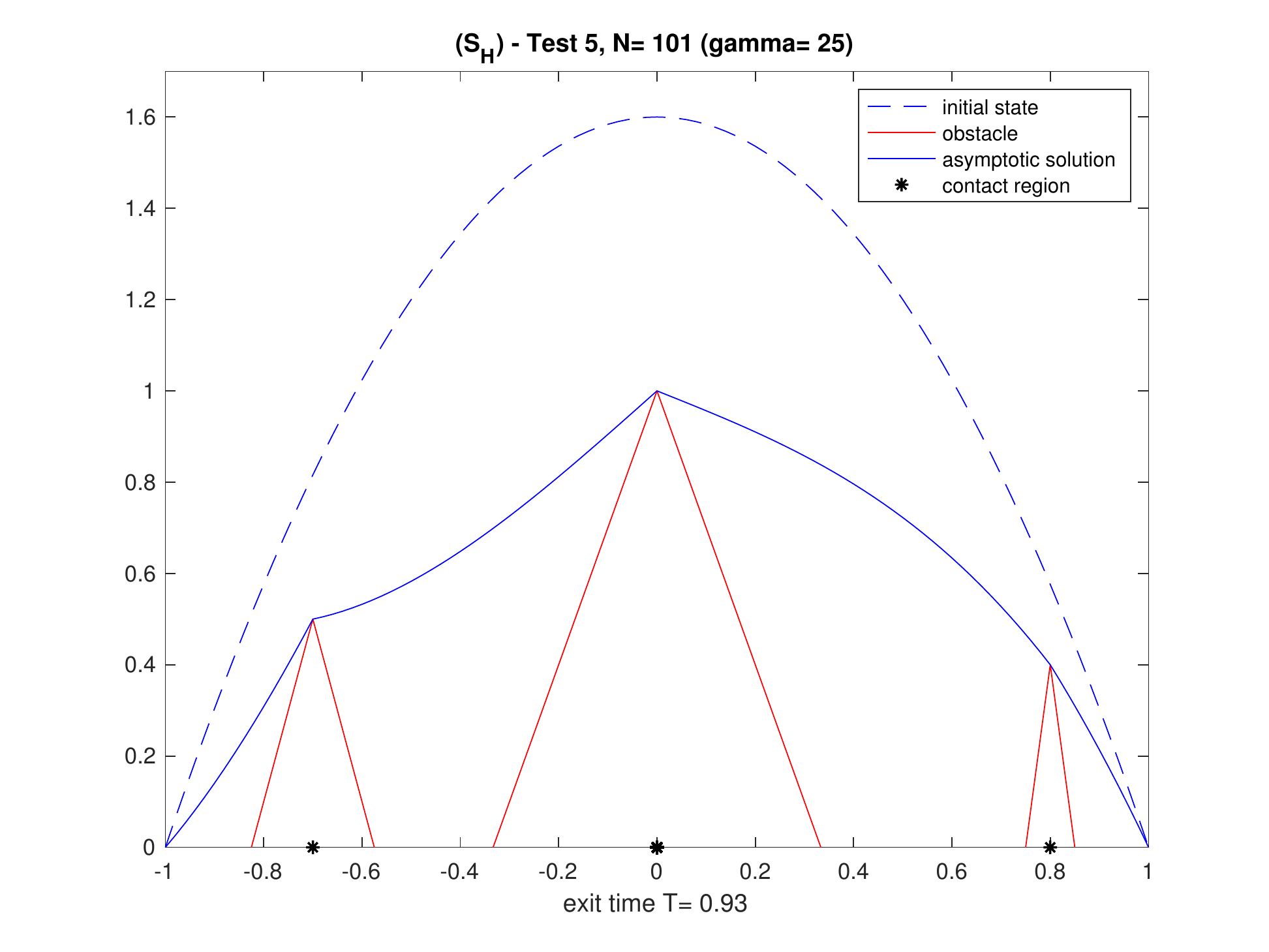}
	\includegraphics[width=4.5cm]{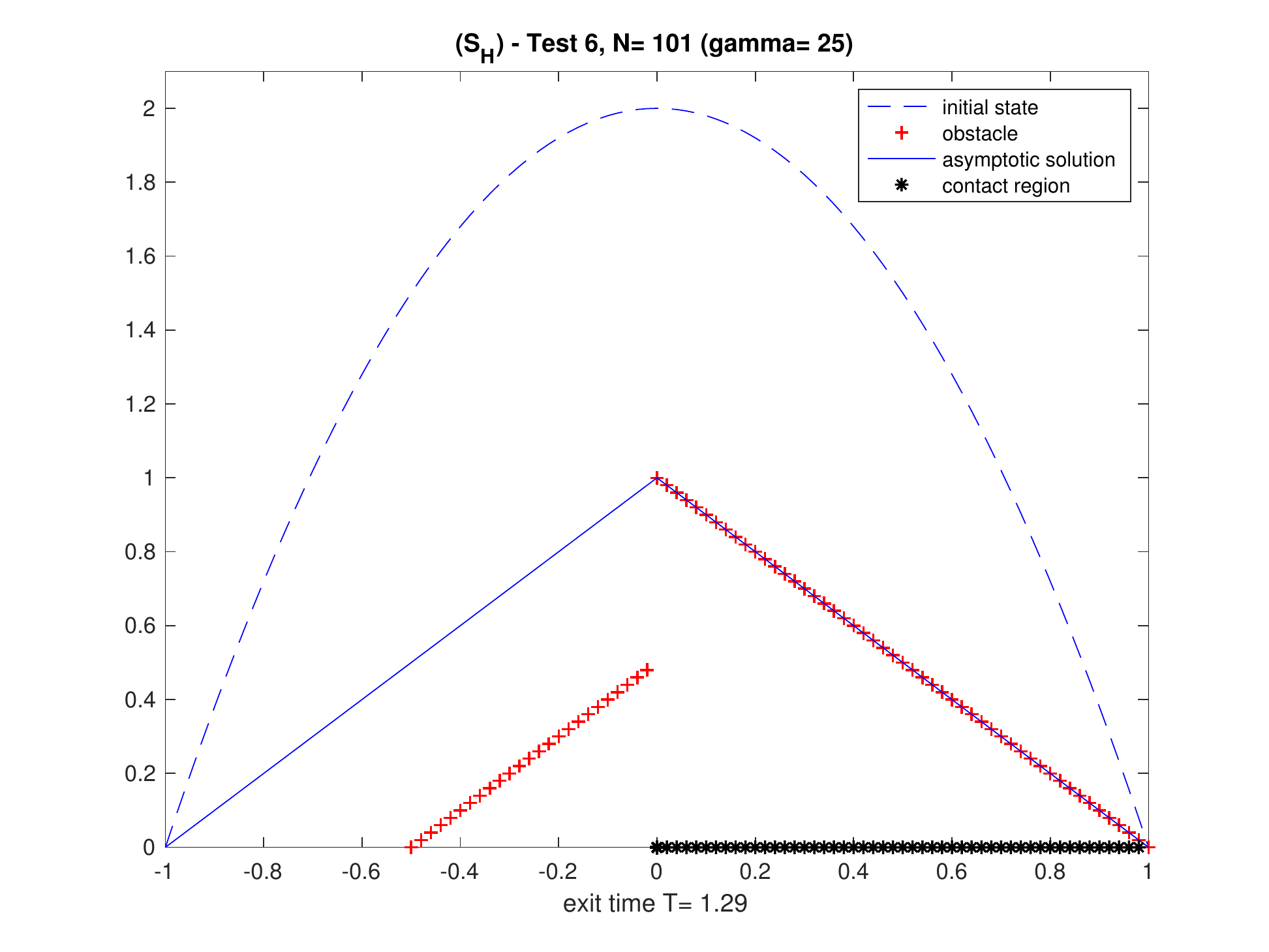}
	\caption{\footnotesize{a) Test 5, b) Test 6. }}	
	\label{F6}
\end{figure}

\medskip{\bf Test 5.}  $u^0=1.6-1.6x^2$, $u^c=\max(1-3|x|,0.5-4|x+0.7|,0.4-8|x-0.8|)$ (three peaks), $f=3x$; the contact set consists of three distinct points (Fig.\ref{F6} a).

\medskip{\bf Test 6.}  $u^0=2-2x^2$, $u^c=x+0.5$ for $x<0$, $u^c=1-x$ for $x\geq0$, $f=0$; $C=[0,1]$  (Fig.\ref{F6} b).

\medskip Finally we report the results of some 2D tests in the square region $\Omega=(-1,1)^2$. For any example we show the final situation, with the surface contact evidence, and explicitely (in blue) the contact area, that is the nodes of the mesh where the solution touches the obstacle.

\begin{figure}[!h]
	\centering      
	\includegraphics[width=4.5cm]{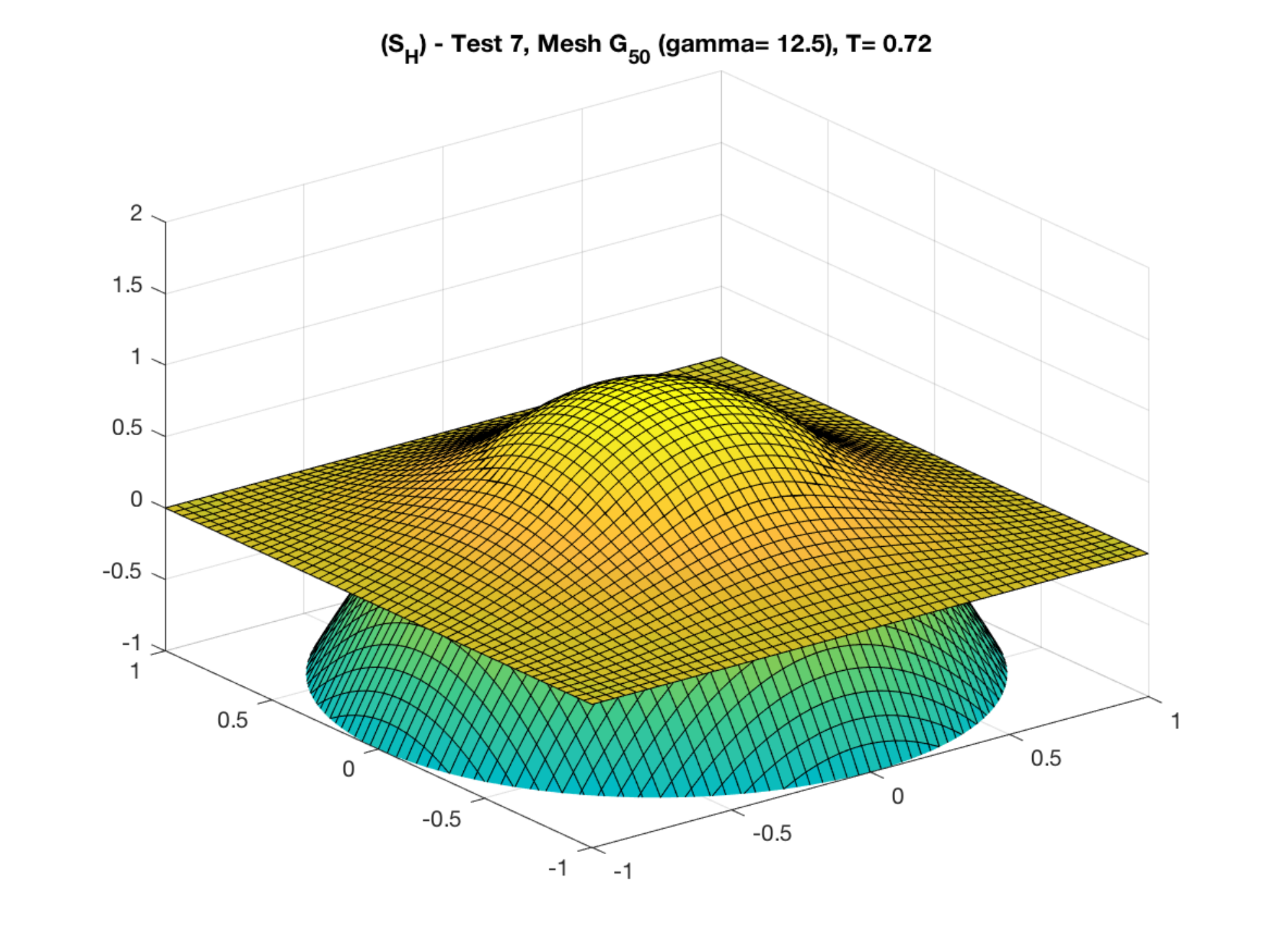}
	\includegraphics[width=4.5cm]{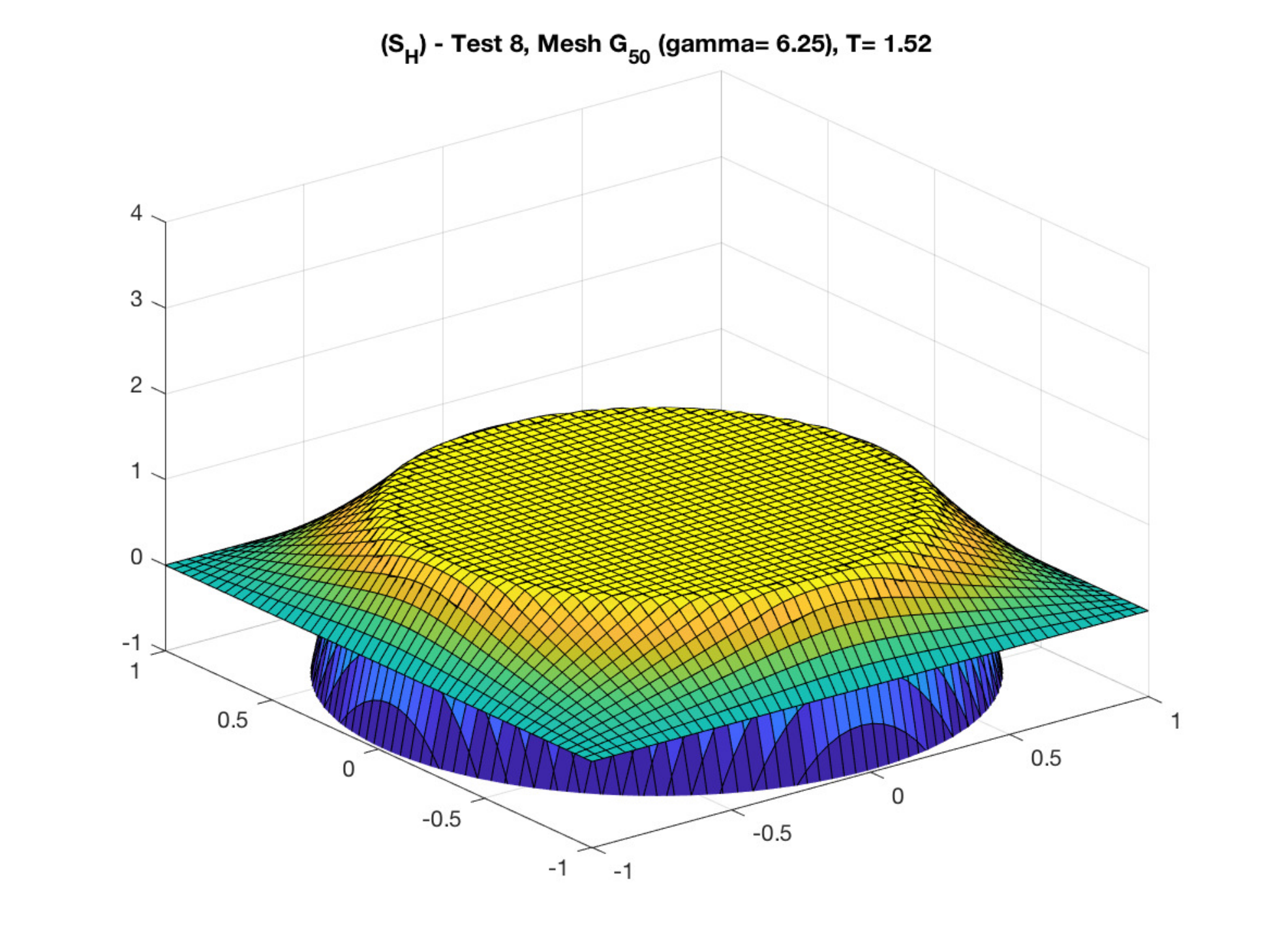}
	\includegraphics[width=4.5cm]{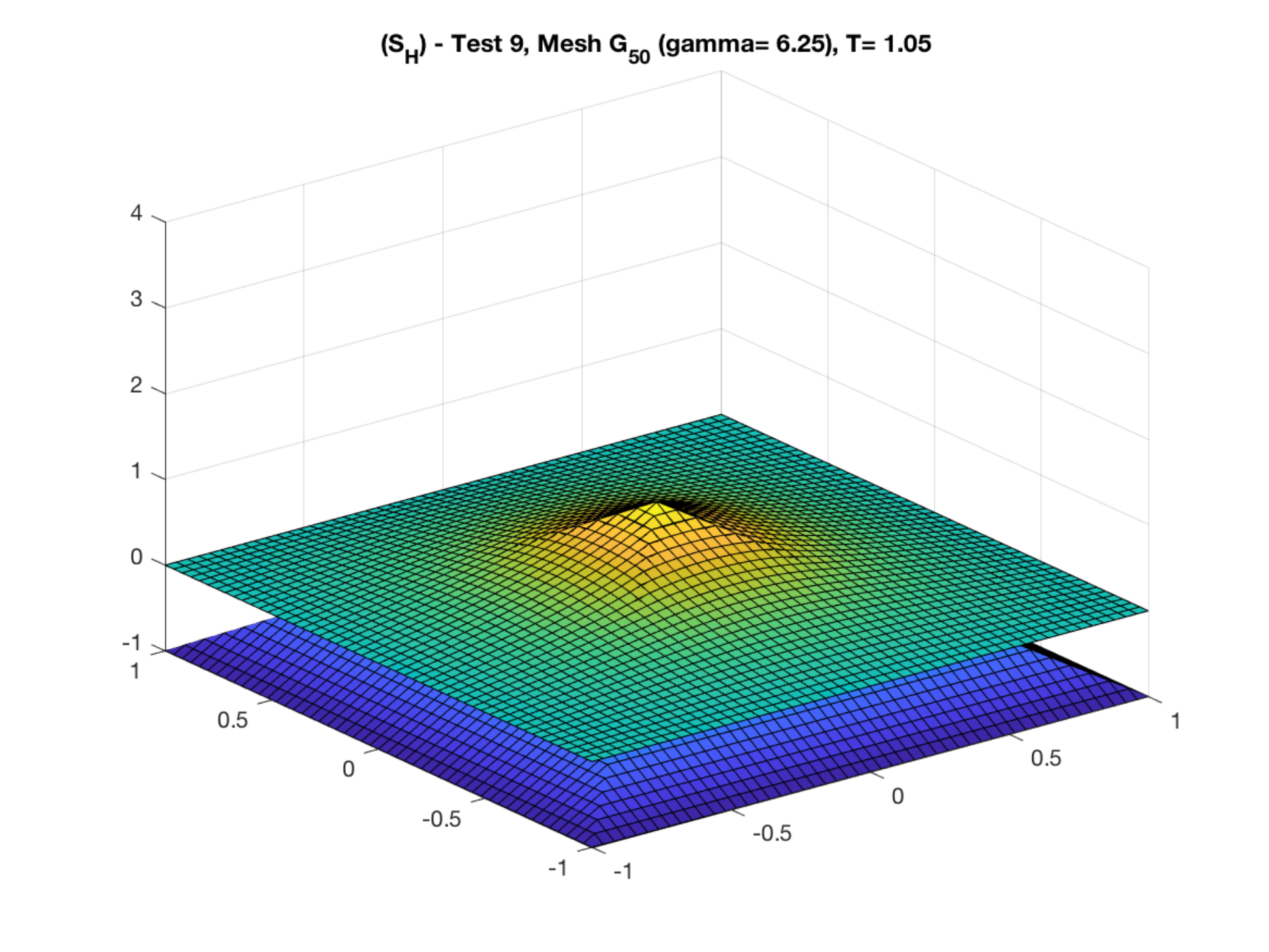}
	\includegraphics[width=4.5cm]{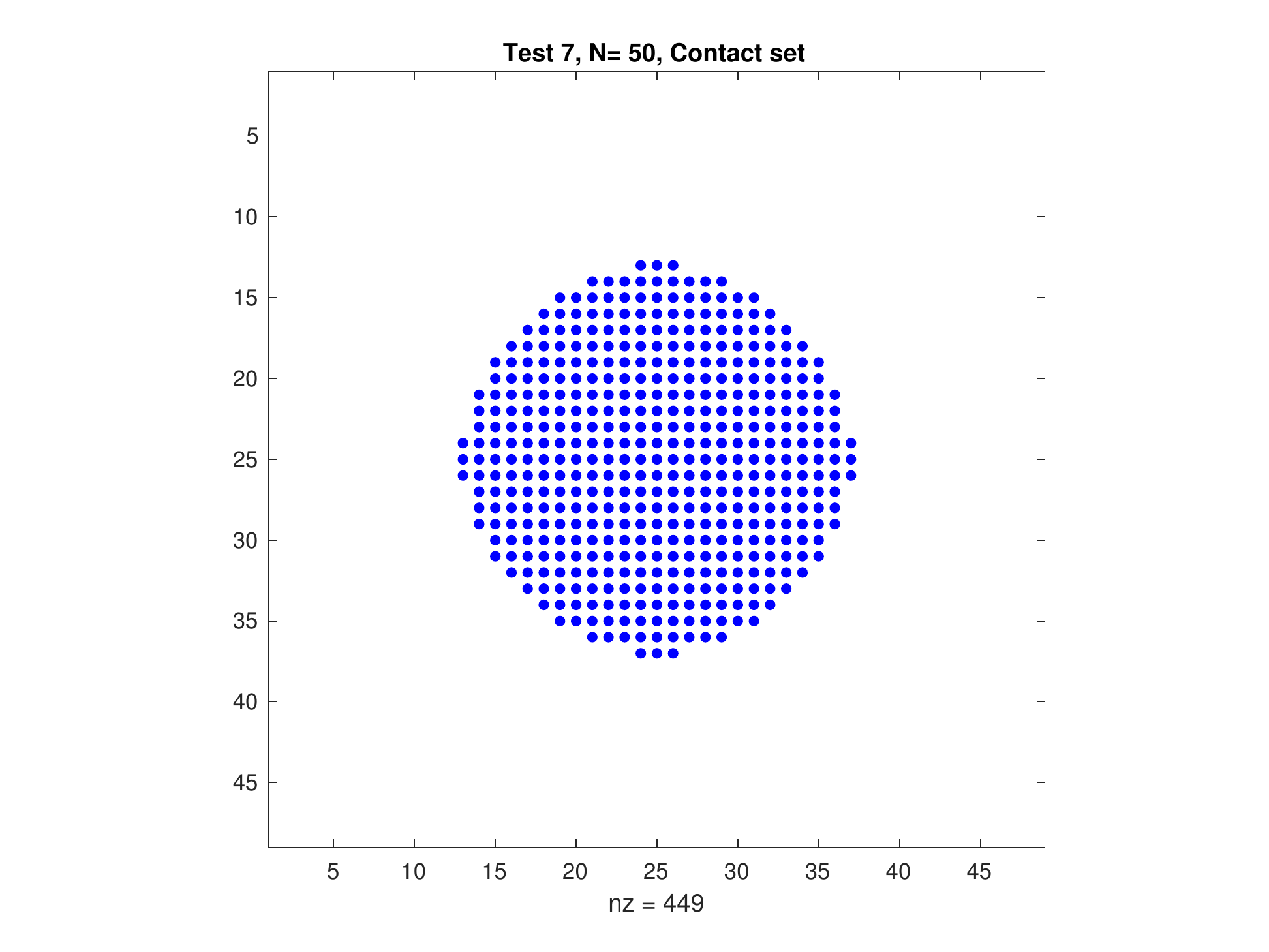}
	\includegraphics[width=4.5cm]{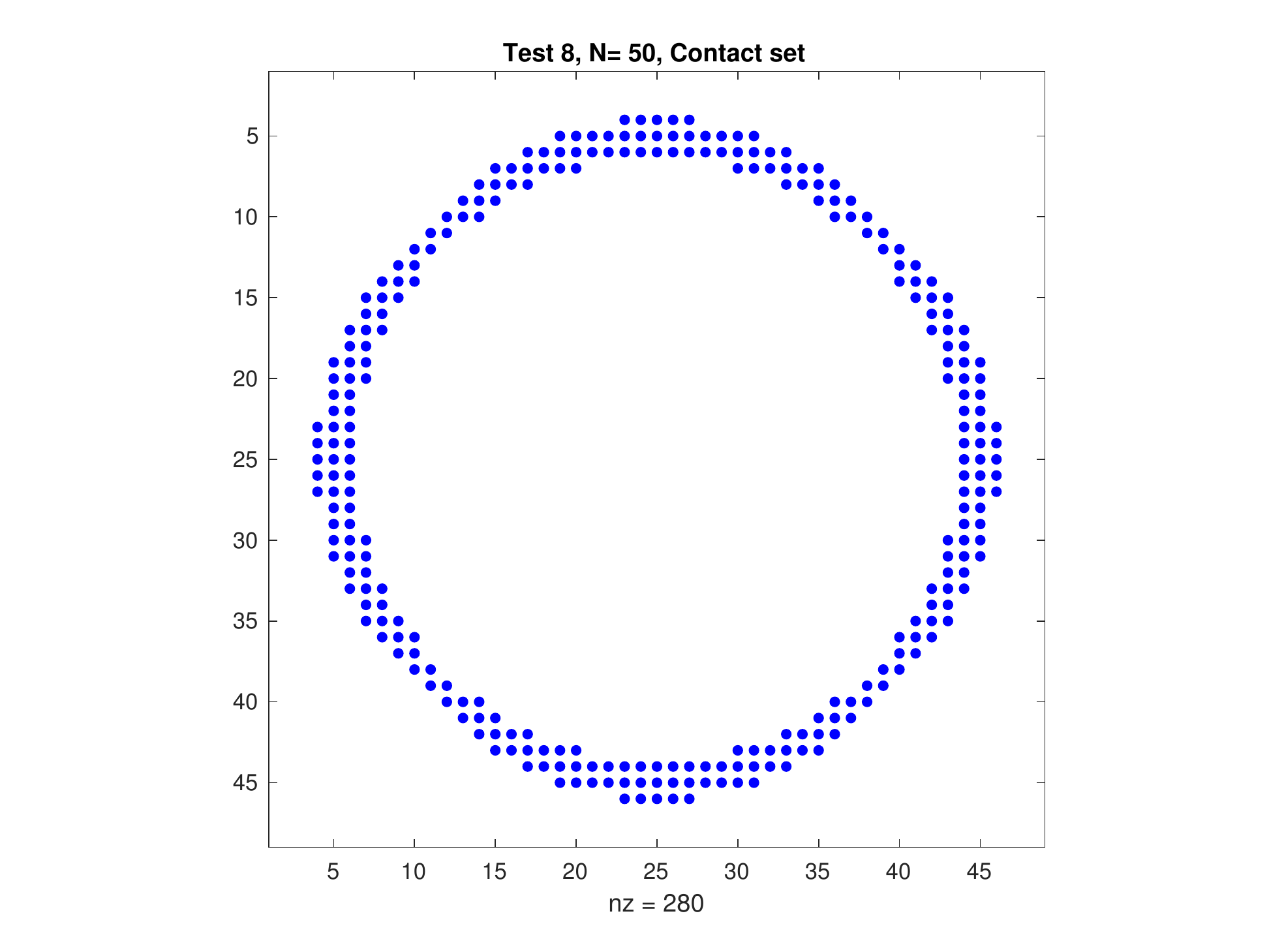}
	\includegraphics[width=4.5cm]{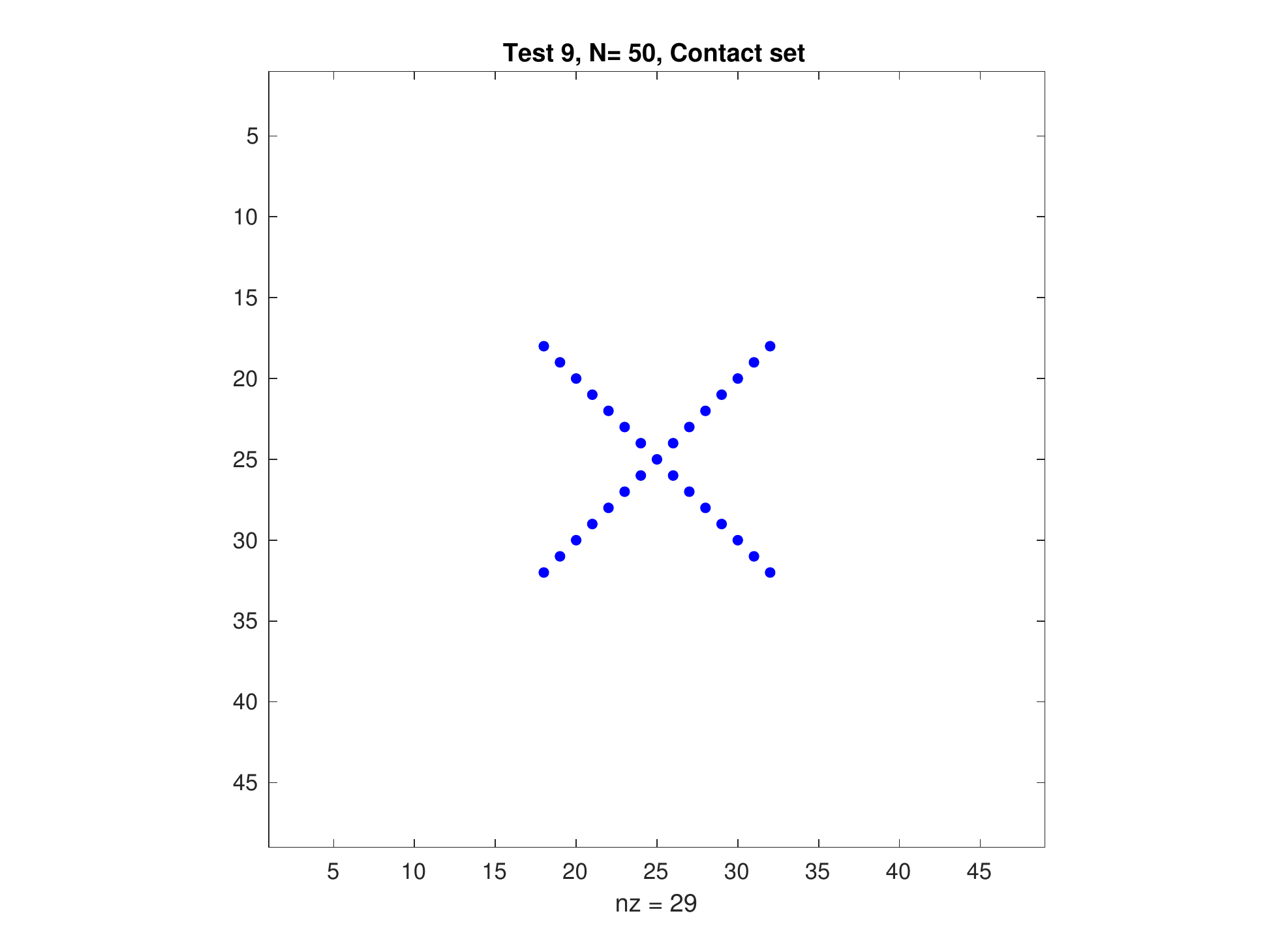}
	\caption{\footnotesize{a) Test 7, b) Test 8, c) Test 9.}}	
	\label{F7}
\end{figure}

\medskip{\bf Test 7.}  $u^0=2(1-x^2)(1-y^2)$;  $u^c=1-2(x^2+y^2)$ (a reversed paraboloid); $f=-1$. The contact set is a disk (Fig.\ref{F7} a).

\medskip{\bf Test 8.} $u^0=4(1-x^2)(1-y^2)$; $u^c=1-(3.5(x^2+y^2)-2)^2$ (a sort of crater of a volcano); $f=0$. The contact set is a circular crown (Fig.\ref{F7} b).

\medskip{\bf Test 9.}  $u^0=2(2-|x+y|-|y-x|)$;  $u^c=(2-|x+y|)-|y-x|)-1$ (a central pyramid); $f=0$. The contact set is made by two crossing lines (Fig.\ref{F7} c).

\begin{figure}[!h]
	\centering      
	\includegraphics[width=4.5cm]{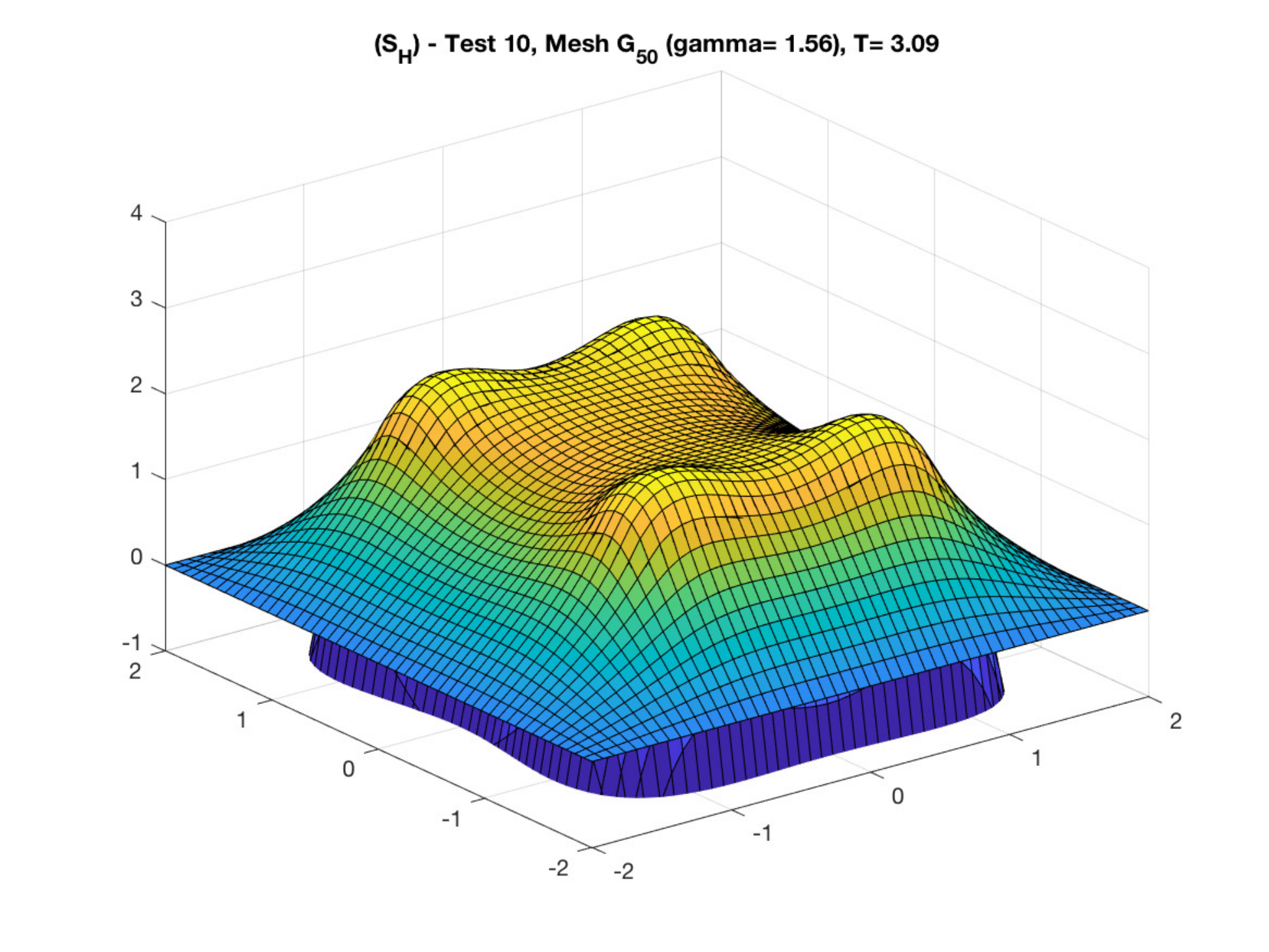}
	\includegraphics[width=4.5cm]{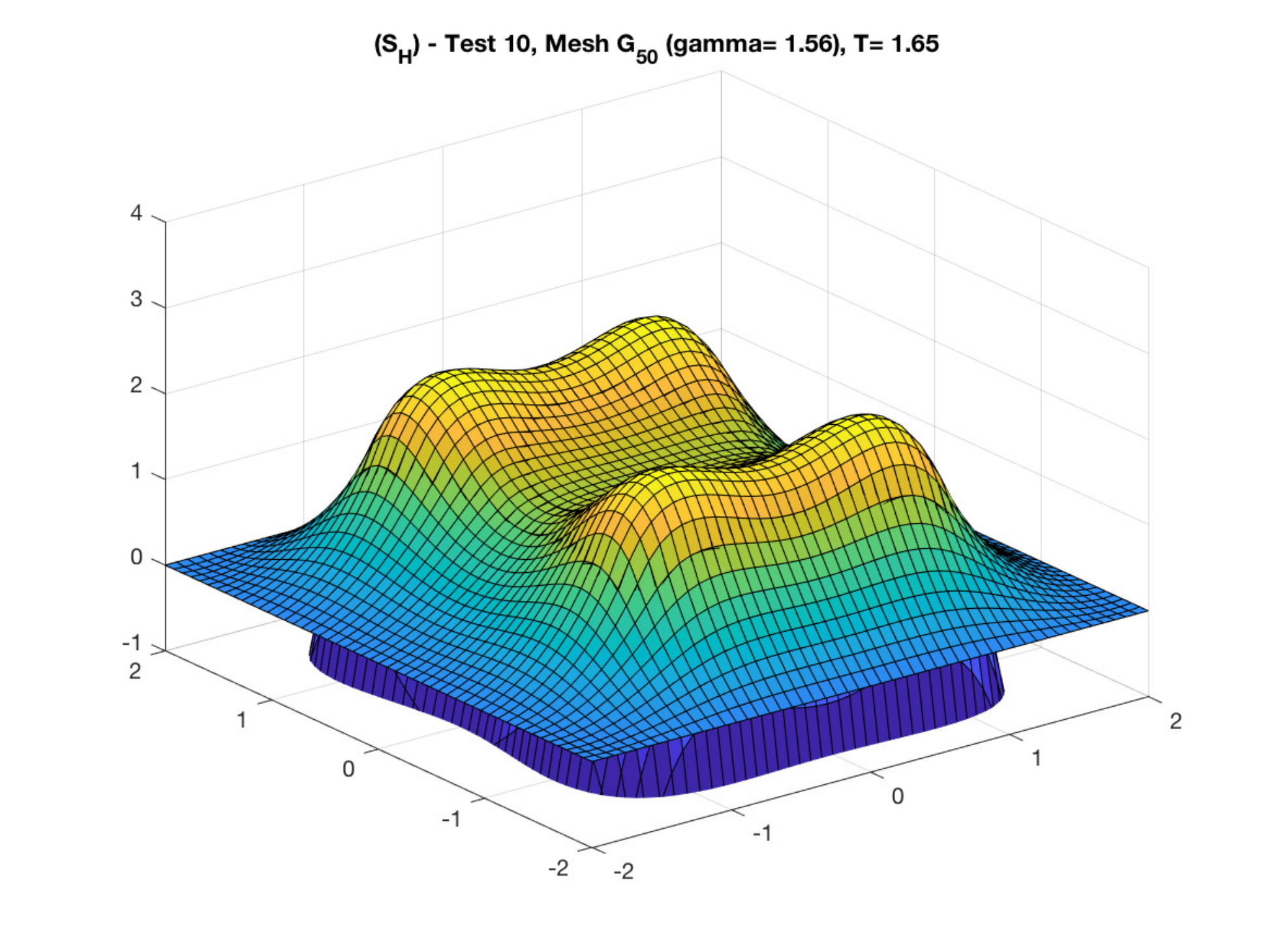}\\
	\includegraphics[width=4.5cm]{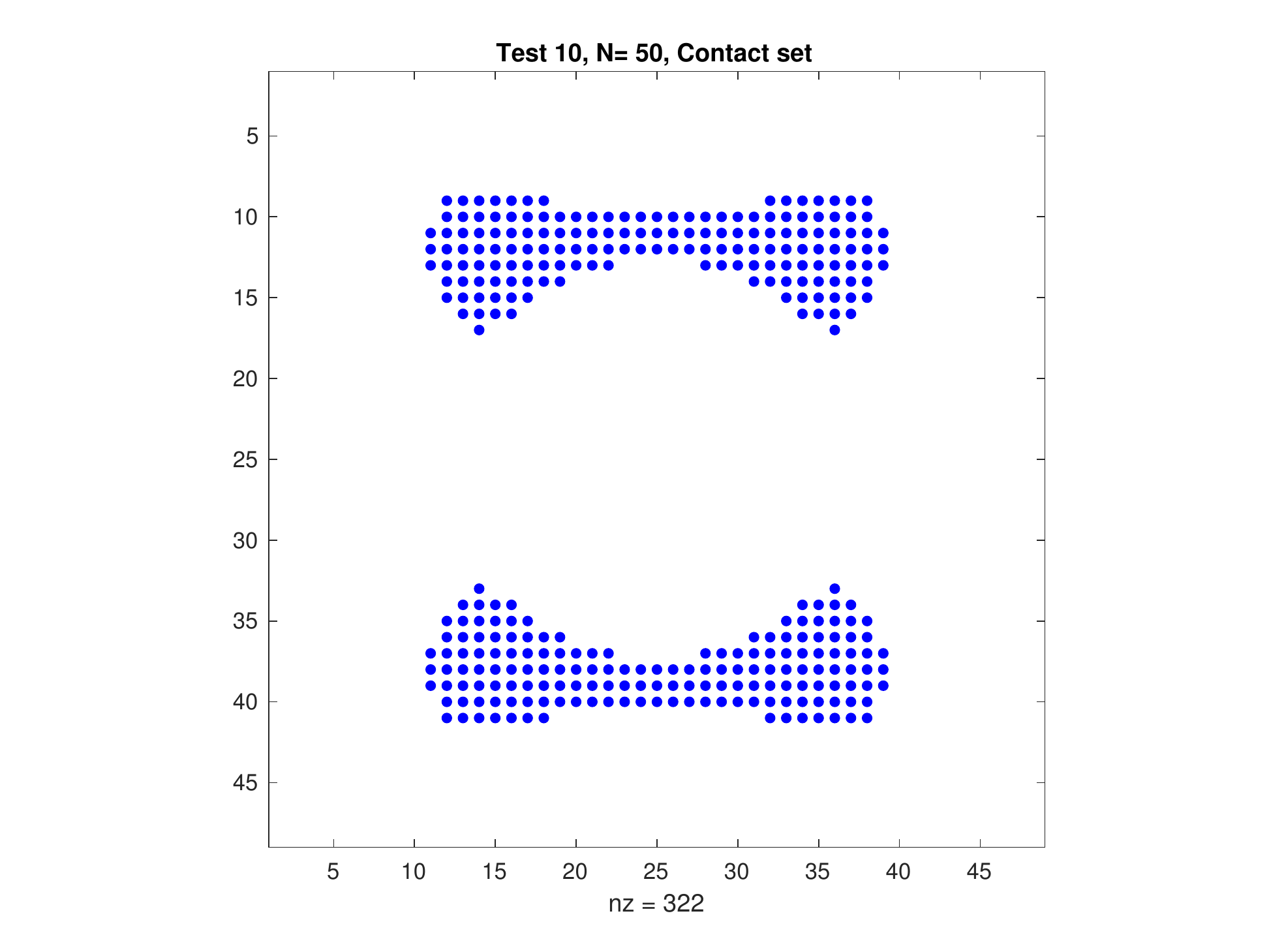}
	\includegraphics[width=4.5cm]{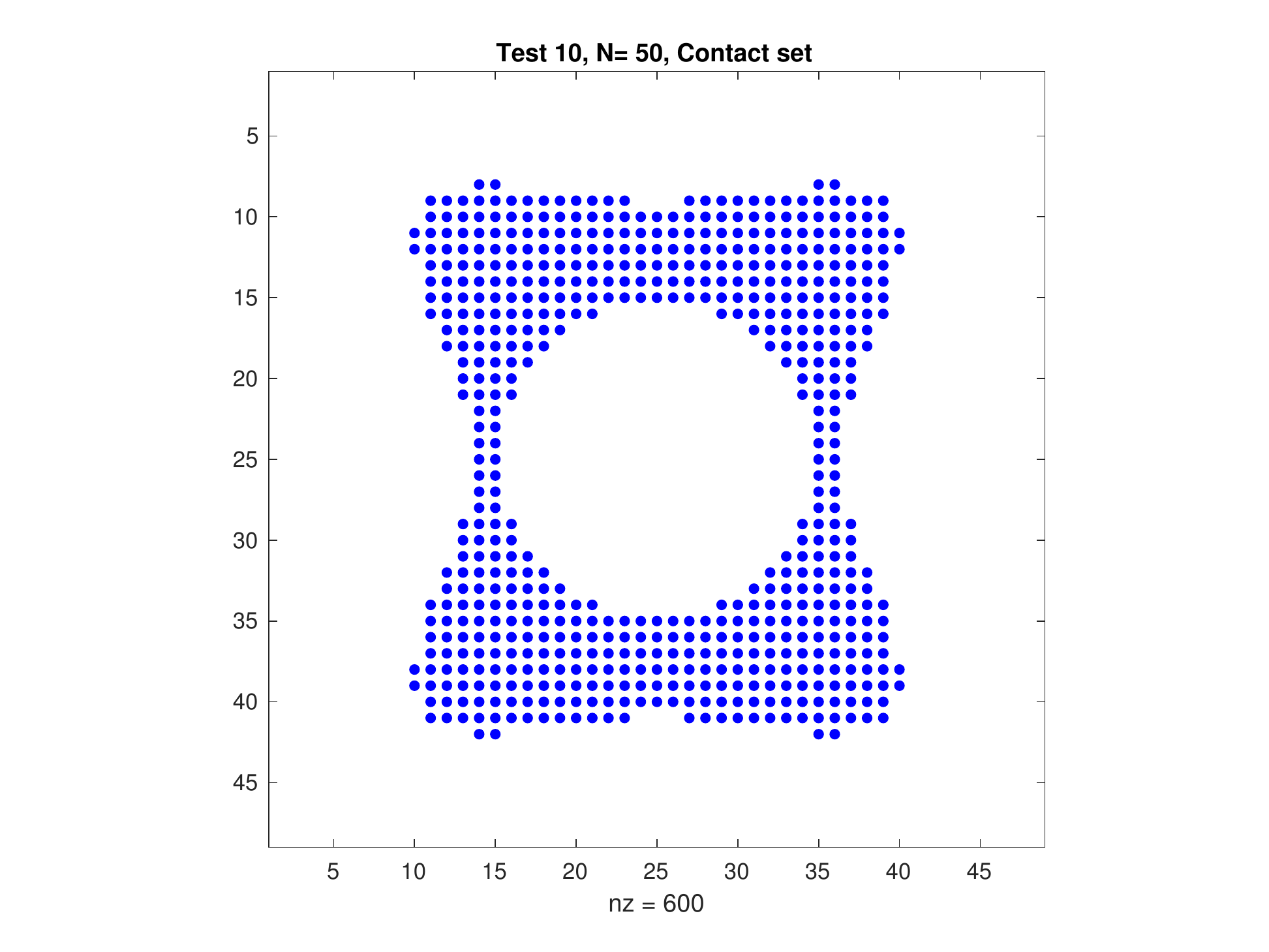}
	\caption{\footnotesize{Test 10. a) $f=0$,  b) $f=-2$.}}	
	\label{F8}
\end{figure}

\medskip{\bf Test 10.}  $u^0=(2-0.5x^2)(2-0.5y^2)$; $u^c=1+x^2+2y^2-x^4-y^4$ (a sort of landscape with hills and valleys);  we compare the final results for $f=0$ and $f=-2$, respectively, with disconnected and connected contact sets (Fig.\ref{F8}). 
%\vfill\eject

\bigskip 

%\textbf{{Acknowledgement.} } Raffaela Capitanelli is member of GNAMPA (INdAM).
%%%%%%%%%%%%%%%%%%%%%%%%%%%%%%%%%%%%

%\newpage

\end{document}